\documentclass[11pt,reqno]{amsart}
\usepackage{fixltx2e}                    
\usepackage[osf]{mathpazo}  
\usepackage[utf8]{inputenc}            
\usepackage{amsmath}                     
\usepackage{amssymb, latexsym, stmaryrd, amsthm, dsfont, amsfonts, amsbsy,amsthm, amsmath, mathrsfs}            
\usepackage{mathtools}                   
\usepackage{bm}                          
\usepackage{enumerate}                   
\usepackage{verbatim}                    
\usepackage{url}   
\usepackage{lscape}                      
\usepackage{soul}

\usepackage[margin=1.3in]{geometry}

\usepackage{microtype}  
\usepackage[all, knot]{xy}
        \xyoption{arc} 
        \xyoption{web}                 
\makeatletter                            
\def\MT@register@subst@font{\MT@exp@one@n\MT@in@clist\font@name\MT@font@list
 \ifMT@inlist@\else\xdef\MT@font@list{\MT@font@list\font@name,}\fi}
\makeatother

\usepackage[pdftex,bookmarks,bookmarksnumbered,linktocpage,   
         colorlinks,linkcolor=blue,citecolor=blue]{hyperref}

\newcommand{\bit}{\begin{itemize}}    
\newcommand{\eit}{\end{itemize}}
\newcommand{\ben}{\begin{enumerate}}
\newcommand{\een}{\end{enumerate}}

\newcommand{\benroman}{\ben[\normalfont (i)]}  
\let\eroman\een

\newcommand{\bde}{\begin{description}}
\newcommand{\ede}{\end{description}}

\newcommand{\?}{\ensuremath{\mkern0.4\thinmuskip}}   
\let\models=\vDash                          

\let\leq=\leqslant
\let\geq=\geqslant

\let\epsilon=\varepsilon
\let\Lambda\varLambda
\let\Gamma\varGamma
\let\Delta\varDelta
\let\Lambda\varLambda
\let\Omega\varOmega
\let\Theta\varTheta
\let\Xi\varXi
\let\Pi\varPi
\let\Sigma\varSigma

\let\clsys=\mathcal                             
\let\class=\mathsf                              
\let\oper=\mathbb                               

\bmdefine{\A}{A}                                
\bmdefine{\2}{2}
\bmdefine{\B}{B}
\bmdefine{\D}{D}
\bmdefine{\M}{M}                                
\bmdefine{\LLL}{L}                              
\bmdefine{\Fm}{Fm}                              
\bmdefine{\zerou}{[0{,}1]}  
\bmdefine{\T}{T}                                

\newcommand{\FF}{\clsys{F}}

\newcommand{\LL}{\mathscr{L}}                   

\newcommand{\PSD}{\oper{P}_{\!\textsc{sd}}^{}}

\newcommand{\SSS}{\oper{S}}

\newcommand{\RRR}{\oper{R}}

\newcommand{\Con}{\mathrm{Con}}                            
\newcommand{\Var}{\mathnormal{V\mkern-.8\thinmuskip ar}} 
\newcommand{\FFi}{\FF\mkern-.5\thinmuskip\mathnormal{i}}  
      
\newcommand{\Mod}{\class{Mod}}

\bmdefine{\boldstar}{\mathchoice{\textstyle*}{\textstyle*}{\textstyle*}{\scriptstyle*}}

\newcommand{\ModS}{\Mod^{\equiv}}

\bmdefine{\btau}{\tau}                                  
\bmdefine{\brho}{\rho}                                  

\bmdefine{\leibniz}{\Omega}        

\bmdefine{\frege}{\Lambda}         

\makeatletter
\newcommand{\tarskidsp}{\mathord%
   {\m@th\raisebox{0pt}[0pt][0pt]{$\stackrel%
   {\raisebox{-2.7pt}[0ex][0pt]{$\displaystyle \,\?\thicksim$}}%
   {\displaystyle\leibniz}$}}}
\newcommand{\tarskitxt}{\mathord%
   {\m@th\raisebox{0pt}[0pt][0pt]{$\stackrel%
   {\raisebox{-2.7pt}[0ex][0pt]{$\,\?\thicksim$}}{\displaystyle\leibniz}$}}}
\newcommand{\tarskiscr}{\mathord%
   {{\m@th\raisebox{0pt}[0pt][0pt]{$\stackrel%
   {\raisebox{-2.4pt}[0ex][0pt]{$\scriptstyle \,\?\thicksim$}}%
   {\scriptstyle\leibniz}$}}}}
\newcommand{\tarskiscrscr}{\mathord%
   {{\m@th\raisebox{0pt}[0pt][0pt]{$\stackrel%
   {\raisebox{-2pt}[0ex][0pt]{$\scriptscriptstyle \,\?\thicksim$}}%
   {\scriptscriptstyle\leibniz}$}}}}
\newcommand{\tarski}{\@ifnextchar ^ %
   {\mathchoice{\tarskidsp\kern-.07em}{\tarskitxt\kern-.07em}%
   {\tarskiscr\kern-.07em}{\tarskiscrscr\kern-.07em}}%
   {\mathchoice{\tarskidsp}{\tarskitxt}{\tarskiscr}{\tarskiscrscr}}}
\makeatother

\theoremstyle{theorem}
\newtheorem{Theorem}{Theorem}[section]
\newtheorem{Proposition}[Theorem]{Proposition}
\newtheorem{Lemma}[Theorem]{Lemma}
\newtheorem{Corollary}[Theorem]{Corollary}

\theoremstyle{definition}
\newtheorem{law}[Theorem]{Definition}
\newtheorem{exa}[Theorem]{Example}
\newtheorem{Fact}{Fact}

\theoremstyle{remark}

\newtheorem{problem}{\bf Problem}
\newtheorem{Remark}[Theorem]{Remark}

\newcommand{\C}{\boldsymbol{C}} 

\begin{document}
\title[The poset of all logics II: Leibniz classes and hierarchy]{The poset of all logics II:\\ Leibniz classes and hierarchy}

\author{R. Jansana and T. Moraschini}
\address{Department of Logic History, and Philosophy of Science, Faculty of Philosophy, University of Barcelona, Carrer de Montalegre $6$, $08001$, Barcelona, Spain}
\email{jansana@ub.edu}
\address{Institute of Computer Science, Academy of Sciences of Czech Republic, Pod Vod\'arenskou v\v{e}\v{z}\'{i} $271/2$, $182$ $07$ Prague $8$, Czech Republic}
\email{moraschini@cs.cas.cz}
\date{\today}

\maketitle

\begin{abstract}
A Leibniz class is a class of logics closed under the formation of term-equivalent logics, compatible expansions, and non-indexed products of sets of logics. We study the complete lattice of all Leibniz classes, called the Leibniz hierarchy. In particular, it is proved that the classes of truth-equational and assertional logics are meet-prime in the Leibniz hierarchy, while the classes of protoalgebraic and equivalential logics are meet-reducible. However, the last two classes are shown to be determined by Leibniz conditions consisting of meet-prime logics only.
\end{abstract}

\section{Introduction}

When ordered under \textit{intepretability} \cite{JaMor19-1}, the class of (propositional) logics forms a preorder. Its associated partial order $\class{Log}$, called the \textit{poset of all logics}, consists of equivalence classes of equi-interpretable logics.\ Building on this formalism, in this paper we introduce and study the notion of a \textit{Leibniz class} of logics. 

From an order-theoretic point of view, Leibniz classes are classes of logics that can be faithfully identified with the upsets of $\class{Log}$ that are closed under infima of arbitrarily large sets. Equivalently, they can be characterized in terms of closure properties as the classes of logics closed under the formation of term-equivalent logics, compatible expansions, and non-indexed products of sets of logics (see \cite{JaMor19-1} for the relevant definitions).

Part of the interest of Leibniz classes lies in the fact that they allow to clarify the yet informal concept of the \textit{Leibniz hierarchy}, i.e.\ a taxonomy in which logics are classified in terms of syntactic principles (up to now recognized on empirical grounds) corresponding to the behaviour of the so-called Leibniz operator \cite{Cz01,AAL-AIT-f,FJaP03b,FJa09}. More precisely, the road we take is to identify the Leibniz hierarchy with the complete lattice of Leibniz classes ordered under inclusion. Remarkably, this abstraction preserves the fact that Leibniz classes are collections of logics globally satisfying some syntactic principles, here called \textit{Leibniz conditions}, consisting of special sequences of logics indexed by all ordinals. 

One of the main advantages of this point of view is that it allowed to unify in \cite{JaMor19-3} the theory of the Leibniz hierarchy with that of the \textit{Maltsev hierarchy} of universal algebra, i.e.\ a taxonomy in which varieties are classified by means of syntactic principles related to the structure of congruence lattices  \cite{HoMcKe88,KeKi06,Neum74,Pixley15can,Tay73,Taylor77,Wille70}.

In this paper we restrict our attention to the study of Leibniz classes. First we show that the majority of classes of logics traditionally associated with the Leibniz hierarchy happen to be Leibniz classes. Among them we count the classes of protoalgebraic, equivalential, and assertional logics, whose respective Leibniz conditions are also found (Theorems \ref{Thm:proto-class}, \ref{Thm:equiv-class}, and \ref{Thm:StrongAssertional}). Then we begin an order-theoretic investigation of the Leibniz hierarchy, understood as the complete lattice of all Leibniz classes. More precisely, we focus on the problem of determining whether the most prominent Leibniz classes are meet-prime or meet-irreducible elements of the Leibniz hierarchy. Affirmative answers to these questions can then be interpreted as stating that the Leibniz classes under consideration capture  primitive or fundamental concepts. Similar problems were studied in the setting of the Maltsev hierarchy for instance in \cite{GaTa84} (see also \cite{BenSeq14,Oprsal18,Tschantz96}).

Some of our results can be summarized as follows. On the one hand the Leibniz classes of truth-equational and assertional logics are shown to be meet-prime (Theorems \ref{Thm:truth-eq-prime} and \ref{Thm:asser-prime}). On the other hand, it is proved that the classes of protoalgebraic and equivalential logics are meet-reducible (Theorems \ref{Thm:dec-equivalential} and \ref{Thm:dec-equivalential}). This negative result is compensated by the observation that the Leibniz condition determining the class of protoalgebraic (resp.\ equivalential) logics consists of logics whose equivalence classes are meet-prime in the poset $\class{Log}$ (Theorems \ref{Thm:basic-proto-prime} and \ref{Thm:basic-equiv-prime}).

\section{Leibniz conditions}

We use the same notation as in \cite{JaMor19-1}. Recall that $\class{OR}$ is the class of ordinals.

\begin{law}
A \textit{strong Leibniz condition} $\Phi$ is simply a logic $\vdash_{\Phi}$.\ A logic $\vdash$ is said to \textit{satisfy} $\Phi$ if ${\vdash_{\Phi}} \leq {\vdash}$, and the class of logics satisfying $\Phi$ is denoted by $\mathrm{Log}(\Phi)$. Similarly, a  \textit{Leibniz condition} $\Phi$ is a class $\{\Phi_\alpha \colon \alpha \in \class{OR}\}$ of strong Leibniz conditions such that if $\alpha \leq \beta$, then $\vdash_{\Phi_\beta} \leq {\vdash_{\Phi_\alpha}}$.\ A logic $\vdash$ is said to \textit{satisfy} $\Phi$ if ${\vdash_{\Phi_{\alpha}}} \leq {\vdash}$ for some $\alpha \in \class{OR}$, and the class of logics satisfying $\Phi$ is denoted by $\mathrm{Log}(\Phi)$.

Accordingly, a class $\oper{K}$ of logics is a (resp.\ \textit{strong}) \textit{Leibniz class} if it is of the form $\mathrm{Log}(\Phi)$ for some (resp.\ strong) Leibniz condition $\Phi$.
\qed
\end{law}

Recall that $	\class{Log}$ is a set-complete meet-semilattice \cite[Thm.\ 4.6]{JaMor19-1}. A subcollection $F \subseteq \class{Log}$ is a \textit{set-complete filter} if it is a non-empty upset closed under infima of sets. Similarly, $F \subseteq \class{Log}$ is a \textit{principal filter} if it is a principal upset of $\class{Log}$. Finally, given a class $\oper{K}$ of logics, we set
 \[
\oper{K}^{\dagger} \coloneqq \{ \llbracket \vdash \rrbracket \colon {\vdash} \in \oper{K} \}.
\]
The following result  is instrumental to construct concrete Leibniz classes.

\begin{Theorem}\label{Thm:Main}
Let $\oper{K}$ be a class of logics. The following conditions are equivalent:
\benroman
\item $\oper{K}$ is a Leibniz class.
\item $\oper{K}$ is closed under term-equivalence, compatible expansions, and non-indexed products of sets of logics.
\item The collection $\oper{K}^{\dagger}$ is a set-complete filter of $\class{Log}$, and $\oper{K} = \{ \vdash \colon \llbracket \vdash \rrbracket \in \oper{K}^{\dagger}\}$.
\eroman
\end{Theorem}
\begin{proof}
We rely on the observation that, given a family $\{ \vdash_{i} \colon i \in I \}$ of logics, the infimum of $\{ \llbracket \vdash_{i} \rrbracket \colon  i \in I \}$ in $\class{Log}$ is $\llbracket \bigotimes_{i \in I}{\vdash_{i}}\rrbracket$ \cite[Thm.\ 4.6]{JaMor19-1}.

(i)$\Rightarrow$(iii): Suppose that $\oper{K}$ is a Leibniz class, i.e.\ that there is a Leibniz condition $\Phi = \{\Phi_{\alpha} \colon \alpha \in \class{OR}\}$ such that $\oper{K} = \class{Log}(\Phi)$. We begin by proving that the collection $\oper{K}^{\dagger}$ is a set-complete filter of $\class{Log}$. 

It is clear that $\oper{K}^{\dagger}$ is a non-empty upset of $\class{Log}$. Then consider a set of logics $\{\vdash_i \colon i \in I \} \subseteq \oper{K}$. Since $\oper{K} = \mathrm{Log}(\Phi)$, for every $i \in I$ there exists an ordinal  $\alpha_i$ such that ${\vdash_{\Phi_{\alpha_i}}} \leq {\vdash_i}$. Let $\beta$ be the supremum of $\{\alpha_i \colon i \in I\}$. Since $\Phi$ is a Leibniz condition, we have ${\vdash_{\Phi_{\beta}}} \leq {\vdash_{\Phi_{\alpha_{i}}}} \leq {\vdash_{i}}$ for every $i \in I$. Hence we obtain ${\vdash_{\Phi_{\beta}}} \leq {\bigotimes_{i \in I}{\vdash_{i}}}$ which, together with $\oper{K} = \mathrm{Log}(\Phi)$, implies that ${\bigotimes_{i \in I}{\vdash_{i}}} \in \oper{K}$. Since  $\llbracket \bigotimes_{i \in I}{\vdash_i} \rrbracket$ is the infimum of $\{ \llbracket \vdash_i \rrbracket \colon i \in I \}$ in $\class{Log}$, we conclude that $\oper{K}^{\dagger}$ is closed under infima of sets. This establishes that $\oper{K}^{\dagger}$ is a set-complete filter of $\class{Log}$.

Now, from the definition of $\oper{K}^{\dagger}$ it follows that $\oper{K} \subseteq \{ \vdash \colon \llbracket \vdash \rrbracket \in \oper{K}^{\dagger}\}$. To prove the other inclusion, consider a logic $\vdash$ such that $\llbracket \vdash \rrbracket \in \oper{K}^{\dagger}$. By the definition of $\oper{K}^{\dagger}$, there is a logic ${\vdash'} \in \oper{K}$ such that ${\vdash'} \leq {\vdash}$. Since ${\vdash'} \in \oper{K}$, there is an ordinal $\alpha$ such that ${\vdash_{\Phi_{\alpha}}} \leq {\vdash'} \leq {\vdash}$. Since $\oper{K} = \mathrm{Log}(\Phi)$, this implies that ${\vdash} \in \oper{K}$ as desired.

(iii)$\Rightarrow$(ii): Suppose that ${\vdash} \in \oper{K}$, and consider a logic $\vdash'$ that is either term-equivalent to $\vdash$ or a compatible expansion of $\vdash$. By \cite[Prop.\ 3.8]{JaMor19-1} we have that ${\vdash} \leq {\vdash'}$ and, therefore, $\llbracket \vdash \rrbracket \leq \llbracket \vdash' \rrbracket$. Since $\oper{K}^{\dagger}$ is an upset of $\class{Log}$ and $\llbracket \vdash \rrbracket \in \oper{K}^{\dagger}$, we obtain that $\llbracket \vdash' \rrbracket \in \oper{K}^{\dagger}$. Together with the fact that $\oper{K} = \{ \vdash'' \colon \llbracket \vdash'' \rrbracket \in \oper{K}^{\dagger}\}$, this yields ${\vdash'} \in \oper{K}$. Hence we conclude that $\oper{K}$ is closed under term-equivalence and compatible expansions.

Then consider a family $\{ \vdash_i \colon i \in I \} \subseteq \oper{K}$. By assumption the infimum of $\{ \llbracket \vdash_{i} \rrbracket \colon i \in I \}$ in $\class{Log}$ belongs to $\oper{K}^{\dagger}$. This amounts to the fact that $\llbracket \bigotimes_{i \in I}{\vdash_{i}}\rrbracket \in \oper{K}^{\dagger}$. As $\oper{K} = \{ \vdash \colon \llbracket \vdash \rrbracket \in \oper{K}^{\dagger}\}$, we conclude that $\bigotimes_{i \in I}{\vdash_{i}} \in \oper{K}$.

(ii)$\Rightarrow$(i): Consider the cumulative hierarchy $\{ V_{\alpha} \colon \alpha \in \class{OR} \}$ of set theory. For every ordinal $\alpha$ we set
\[
\oper{K}_{\alpha} \coloneqq \oper{K} \cap V_{\alpha}\text{ and }\vdash_{\alpha} \coloneqq \bigotimes \oper{K}_{\alpha}.
\]
Note that if $\oper{K}_{\alpha} = \emptyset$, then $\vdash_{\alpha}$ is the inconsistent logic in the empty language over $\omega$ variables.  Also note that if $\alpha \leq \beta$, then $\oper{K}_{\alpha} \subseteq \oper{K}_{\beta}$ and, therefore, ${\vdash_{\beta}} \leq {\vdash_{\alpha}}$. 
In particular, this implies that the following is a Leibniz condition:
\[
\Phi \coloneqq  \{\vdash_{\alpha} \colon \alpha \in \class{OR}\}.
\]

To conclude the proof, it suffices to show that $\mathrm{Log}(\Phi) = \oper{K}$. To prove the right-to-left inclusion, consider a logic ${\vdash} \in \oper{K}$. Since ${\vdash}$ is a set, there is an ordinal $\alpha$ such that ${\vdash} \in V_{\alpha}$. This implies that ${\vdash} \in \oper{K}_{\alpha}$ and, therefore, that ${\vdash_{\alpha}} \leq {\vdash}$. Hence we conclude that ${\vdash} \in \mathrm{Log}(\Phi)$.

To prove the other inclusion, consider ${\vdash} \in \class{Log}(\Phi)$. There exists an ordinal $\alpha$ such that ${\vdash_{\alpha}} \leq {\vdash}$. By \cite[Prop.\ 3.8]{JaMor19-1} this implies that $\vdash$ is term-equivalent to a compatible expansion of $\bigotimes \oper{K}_{\alpha}$. As $\oper{K}_{\alpha} \subseteq \oper{K}$ and $\oper{K}$ is closed under non-indexed products of sets of logics, compatible expansions, and term-equivalence, this implies ${\vdash} \in \oper{K}$.
\end{proof}

\begin{Corollary}\label{Cor:StrongCharac}
A class $\oper{K}$ of logics is a strong Leibniz class if and only if  $\oper{K}^{\dagger}$ is a principal filter of $\class{Log}$, and $\oper{K} = \{ \vdash \colon \llbracket \vdash \rrbracket \in \oper{K}^{\dagger}\}$.
\end{Corollary}

\begin{proof}[Proof sketch.]
The proof of the ``only if'' part is an easier variant of the one of part (i)$\Rightarrow$(iii) of Theorem \ref{Thm:Main}. Then we sketch the ``if'' part only. By the assumption we know that $\oper{K}^{\dagger}$ is the upset generated by $\llbracket \vdash \rrbracket$, for some logic ${\vdash} \in \oper{K}$. Let $\Phi$ be the strong Leibniz condition determined by $\vdash$. It is not hard to see that $\mathrm{Log}(\Phi) = \oper{K}$ and, therefore, that $\oper{K}$ is a strong Leibniz class.
\end{proof}

\begin{Remark}\label{Rem:equivalential}
Typical applications of Theorem \ref{Thm:Main} comprise proofs that certain collections of logics are Leibniz classes. For instance, in \cite[Props.\ 3.8 and 6.1]{JaMor19-1} it was shown that the class of equivalential logics is closed under term-equivalence, compatible expansions, and non-indexed products of sets of logics. By Theorem \ref{Thm:Main}(ii) we conclude that it is a Leibniz class.
\qed
\end{Remark}

\begin{Remark}\label{Rem:extensions}
From Theorem \ref{Thm:Main}(ii) it follows that Leibniz classes are closed under the formation of extensions of logics, as these are special cases of compatible expansions.
\qed
\end{Remark}

In this paper we identify the intuitive concept of the \textit{Leibniz hierarchy} with the poset of all Leibniz classes ordered under the inclusion relation.
\begin{Proposition}\label{Prop:meets-and-joins}
The Leibniz hierarchy is a complete lattice in which infima are intersections and and suprema are obtained as follows for every collection $\{ \oper{K}_{i} \colon i \in I \}$ of Leibniz classes, where $I$ can be a proper class:
\[
\bigvee_{i \in I} \oper{K}_{i} = \{ {\vdash} \colon {\vdash} \text{ is a logic and }{\bigotimes_{j \in J}{\vdash_{j}}} \leq {\vdash} \text{ for some subset }\{ \vdash_{j} \colon j \in J \} \subseteq \bigcup_{i \in I}\oper{K}_{i} \}.
\] 
\end{Proposition}

\begin{proof}
Immediate from Theorem \ref{Thm:Main}.
\end{proof}

\begin{Remark}
The statement of Proposition \ref{Prop:meets-and-joins} presupposes the we can meaningfully speak of very large intersections and unions of classes (of logics), as the collection $\{ \oper{K}_{i} \colon i \in I \}$ is in general a collection of classes. However, for our purposes this problem is immaterial as we will only work with finite joins and meets of Leibniz classes.
\qed
\end{Remark}

\section{Examples of Leibniz classes}

In this section we show that a range of classes of logics, traditionally associated with the yet informal concept of the Leibniz hierarchy in abstract algebraic logic, are indeed Leibniz classes.

A logic $\vdash$ is said to be \textit{protoalgebraic} \cite{BP86,Cz01} if there is a non-empty set\footnote{In the literature the set $\Delta$ is not required to be non-empty. However, this restriction is almost immaterial as, in a fixed language, there is a unique protoalgebraic logic with an empty $\Delta$, namely the almost inconsistent logic \cite[Prop.\ 6.11.5]{AAL-AIT-f}.} $\Delta(x, y, \vec{z})$ of formulas such that for every $\langle \A, F \rangle \in \Mod(\vdash)$ and $a, b \in A$,
\[
\langle a, b \rangle \in \leibniz^{\A}F \Longleftrightarrow \Delta^{\A}(a, b, \vec{c}) \subseteq F \text{, for every }\vec{c} \in A.
\]
In this case, we say that $\Delta(x, y, \vec{z})$ is a set of \textit{congruence formulas with parameters} for $\vdash$.

\begin{Theorem}\label{Thm:Protoalgebraic-characterization}
Let $\vdash$ be a logic. The following conditions are equivalent:
\benroman
\item $\vdash$ is protoalgebraic.
\item $\vdash$ has theorems and $\ModS(\vdash) = \RRR(\Mod(\vdash))$.
\item There is a non-empty set of formulas $\nabla(x, y)$ such that 
\[
\emptyset \vdash \nabla(x, x) \qquad x, \nabla(x, y) \vdash y.
\]
\eroman
In this case, the following is a set congruence formulas with parameters for $\vdash$:
\[
\hat{\nabla}(x, y, \vec{z}) \coloneqq \{ \varphi(\psi(x, \vec{z}), \psi(y, \vec{z})) \colon \varphi(x, y) \in \nabla \text{ and }\psi(x, \vec{z}) \in Fm_{\LL_{\vdash}}(\omega) \}.
\]
\end{Theorem}

\begin{proof}
For the equivalence between (i), (ii), and (iii), see \cite[Thms.\ 6.7, 6.17, and 6.57]{AAL-AIT-f}. The fact that $\hat{\nabla}$ is a set of equivalence formulas with parameters for $\vdash$ follows from \cite[Thm.\ 13.5]{BP92} (see also \cite[Prop.\ 3.2]{Font13}).
\end{proof}

Our aim is to prove that protoalgebraic logics form a Leibniz class. To this end, it is convenient to introduce the following concept:

\begin{law}\label{Def:proto-basic}
Given an infinite cardinal $\kappa$, let $\mathscr{L}^{\kappa}_{\class{P}}$ be the language consisting of the binary symbols $\{ \multimap_{\alpha} \colon \alpha < \kappa \}$ and the $n$-ary symbols $\{ \ast_{n \alpha} \colon \alpha < \kappa \}$ for $0 < n \in \omega$. We set
\[
\nabla_{\kappa}(x, y) \coloneqq \{ x \multimap_{\alpha} y \colon \alpha < \kappa \}.
\]
The \textit{basic protoalgebraic logic of rank $\kappa$} is the logic $\vdash_{\class{P}}^{\kappa}$ formulated on $Fm_{\LL_{\class{P}}^{\kappa}}(\omega)$ determined by the rules
\[
\pushQED{\qed} \emptyset \rhd \nabla_{\kappa}(x, x) \qquad x, \nabla_{\kappa}(x, y) \rhd y.\qedhere \popQED
\]
\end{law}

The following result explains the role of $\vdash_{\class{P}}^{\kappa}$.

\begin{Proposition}\label{Prop:proto-class1}  Let $\vdash$ be a logic.
\benroman
\item $\vdash$ is protoalgebraic if and only if ${\vdash_{\class{P}}^{\kappa}} \leq {\vdash}$ for every (equiv.\ some) infinite cardinal $\kappa \geq \vert \LL_{\vdash} \vert$.
\item If $\kappa$ is an infinite cardinal and ${\vdash_{\class{P}}^{\kappa}} \leq {\vdash}$, then $\vdash$ is protoalgebraic.
\eroman
\end{Proposition}

\begin{proof}
(ii): Let $\btau$ be an interpretation of $\vdash_{\class{P}}^{\kappa}$ into $\vdash$. From \cite[Prop.\ 3.3]{JaMor19-1} we obtain $\emptyset \vdash \btau[\nabla_{\kappa}(x, x)]$ and $x, \btau[\nabla_{\kappa}(x, y)] \vdash y$. Then the set $\btau[\nabla_{\kappa}(x, y)]$ witnesses the validity of condition (iii) of Theorem \ref{Thm:Protoalgebraic-characterization}, whence $\vdash$ is protoalgebraic.

(i): As a special instance of (ii) we obtain that if ${\vdash_{\class{P}}^{\kappa}} \leq {\vdash}$ for some infinite cardinal $\kappa \geq \vert \LL_{\vdash} \vert$, then $\vdash$ is protoalgebraic. Then suppose that $\vdash$ is protoalgebraic, and consider any infinite cardinal $\kappa \geq \vert \LL_{\vdash} \vert$. By Theorem \ref{Thm:Protoalgebraic-characterization} there is a set of formulas $\nabla(x, y) \subseteq Fm(\vdash)$ such that
\begin{equation}\label{Eq:Pkappa1}
\emptyset\vdash \nabla(x, x) \text{ and }x, \nabla(x, y) \vdash y.
\end{equation}

Now, observe that $\vert \nabla \vert \leq \max\{\omega, \vert \LL_{\vdash}\vert\} \leq \kappa$. Since $\nabla \ne \emptyset$, there is a surjective map
\[
f \colon \{  \multimap_{\alpha} \colon \alpha < \kappa \} \to \nabla(x_{1}, x_{2}).
\]
Similarly, since $\vert Fm_{\LL_{\vdash}}(\omega) \vert \leq \max\{ \omega, \vert \LL_{\vdash} \vert \}\leq \kappa$, for every $0 < n \in \omega$ there is a surjective map
\[
g_{n} \colon \{ \ast_{n \alpha} \colon \alpha < \kappa \} \to \{ \varphi \in Fm_{\LL_{\vdash}}(\omega) \colon \varphi = \varphi(x_{1}, \dots, x_{n})\}.
\]
Observe that the maps $f$ and $\{ g_{n} \colon 0 < n \in \omega \}$ can be turned in the natural way into a single translation $\btau$ of $\LL^{\kappa}_{\class{P}}$ into $\LL_{\vdash}$.

Recall by Theorem \ref{Thm:Protoalgebraic-characterization} that $\hat{\nabla}_{\kappa}$ and $\hat{\nabla}$ are sets of congruence formulas with parameters for $\vdash_{\class{P}}^{\kappa}$ and $\vdash$ respectively. We claim that 
\begin{equation}\label{Eq:proto-claim-kappa}
\hat{\nabla}(x, y, \vec{z}) = \btau[\hat{\nabla}_{\kappa}(x, y, \vec{z})].
\end{equation}
We begin by proving the inclusion from left to right. Consider $\varphi \in \hat{\nabla}$. There are $\psi(x, y) \in \nabla(x, y)$ and $\gamma(x, z_{1}, \dots, z_{n-1}) \in Fm_{\LL_{\vdash}}(\omega)$ such that
\[
\varphi = \psi(\gamma(x, z_{1}, \dots, z_{n-1}), \gamma(y, z_{1}, \dots, z_{n-1})).
\]
Since $f$ and $g_{n}$ are surjective, there are $\alpha, \beta < \kappa$ such that $\btau(\multimap_{\alpha}) = \psi(x_{1}, x_{2})$ and $\btau(\ast_{n\beta}) = \gamma(x_{1}, \dots, x_{n})$. Clearly,
\[
\varphi = \btau(\ast_{n\beta}(x, z_{1}, \dots, z_{n-1})\multimap_{\alpha}\ast_{n\beta}(y, z_{1}, \dots, z_{n-1}) ).
\]
Moreover, the definition of $\hat{\nabla}_{\kappa}$ guarantees
\[
\ast_{n\beta}(x, z_{1}, \dots, z_{n-1})\multimap_{\alpha}\ast_{n\beta}(y, z_{1}, \dots, z_{n-1}) \in \hat{\nabla}_{\kappa}.
\]
Hence we conclude that $\varphi \in \btau[\hat{\nabla}_{\kappa}]$. This establishes the inclusion from left to right in (\ref{Eq:proto-claim-kappa}).

To prove the other inclusion, consider $\varphi \in \btau[\hat{\nabla}_{\kappa}]$. Then there are are $\alpha < \kappa$ and $\gamma(x, \vec{z}) \in Fm(\vdash_{\class{P}}^{\kappa})$ such that $\varphi = \btau(\gamma(x, \vec{z}) \multimap_{\alpha} \gamma(y, \vec{z}))$. Set
\[
\gamma'(x, \vec{z}) \coloneqq \btau(\gamma(x, \vec{z})) \text{ and }\psi(x, y) \coloneqq \btau(x \multimap_{\alpha} y).
\]
From the definition of $\hat{\nabla}$ and the fact that $f(\multimap_{\alpha})  \in \nabla(x_{1}, x_{2})$, it follows
\[
\varphi = \btau(\gamma(x, \vec{z}) \multimap_{\alpha} \gamma(y, \vec{z})) = \psi(\gamma'(x, \vec{z}), \gamma'(y, \vec{z})) \in \hat{\nabla}.
\]
This establishes the equality in (\ref{Eq:proto-claim-kappa}).

To conclude the proof, it suffices to show that $\btau$ is an interpretation of $\vdash_{\class{P}}^{\kappa}$ into $\vdash$. To this end, consider a matrix $\langle \A, F \rangle \in \ModS(\vdash)$. We begin by showing that $\langle \A^{\btau}, F \rangle$ is a model of $\vdash_{\class{P}}^{\kappa}$. Observe that $\nabla(x, y) = \btau[\nabla_{\kappa}(x, y)]$. Together with (\ref{Eq:Pkappa1}) and $\langle \A, F \rangle \in \Mod(\vdash)$, this yields that the matrix $\langle \A, F \rangle$ is a model of the rules $\emptyset\rhd \btau[\nabla_{\kappa}(x, x)]$ and $x, \btau[\nabla_{\kappa}(x, y)] \rhd y$. As a consequence, $\langle \A^{\btau}, F \rangle$ is a model of the rules $\emptyset \rhd \nabla_{\kappa}(x, x)$ and $x, \nabla_{\kappa}(x, y) \rhd y$, whence it is a model of $\vdash_{\class{P}}^{\kappa}$.

Now, for every $a, b \in A$,
\begin{align*}
a = b &\Longleftrightarrow \langle a, b \rangle \in \leibniz^{\A}F\\
& \Longleftrightarrow \hat{\nabla}^{\A}(a, b, \vec{c})\subseteq F\text{ for every }\vec{c}\in A\\
& \Longleftrightarrow\btau[\hat{\nabla}_{\kappa}]^{\A}(a, b, \vec{c})\subseteq F\text{ for every }\vec{c}\in A\\
& \Longleftrightarrow \hat{\nabla}_{\kappa}^{\A^{\btau}}(a, b, \vec{c}) \subseteq F\text{ for every }\vec{c}\in A\\
&\Longleftrightarrow \langle a, b \rangle \in \leibniz^{\A^{\btau}}F.
\end{align*}
The above equivalences are justified as follows: the first is a consequence of the fact that $\langle \A, F \rangle$ is reduced by Theorem \ref{Thm:Protoalgebraic-characterization}(ii), the second follows from the fact that $\langle \A, F \rangle \in \Mod(\vdash)$ and that $\hat{\nabla}$ is a set of congruence formulas with parameters for $\vdash$, the third from (\ref{Eq:proto-claim-kappa}), the fourth is straightforward, and the fifth from the fact that $\langle \A^{\btau}, F \rangle \in \Mod(\vdash_{\class{P}}^{\kappa})$ and that $\hat{\nabla}_{\kappa}$ is a set of congruence formulas with parameters for $\vdash_{\class{P}}^{\kappa}$. The above display implies that the congruence $\leibniz^{\A^{\btau}}F$ is the identity relation. As a consequence, we obtain $\langle \A, F \rangle \in \RRR(\Mod(\vdash_{\class{P}}^{\kappa})) = \ModS(\vdash_{\class{P}}^{\kappa})$. This establishes that $\btau$ is an interpretation of $\vdash_{\class{P}}^{\kappa}$ into $\vdash$, as desired.
\end{proof}

For every ordinal $\alpha$, let $\vdash_{\class{P}}^{\alpha}$ be the logic $\vdash_{\class{P}}^{\omega + \vert \alpha \vert}$.

\begin{Theorem}\label{Thm:proto-class}
The sequence $\Phi = \{ \vdash_{\class{P}}^{\alpha}  \colon \alpha \in \class{OR} \}$ is a Leibniz condition and $\mathrm{Log}(\Phi)$ coincides with the class of protoalgebraic logics. In particular, protoalgebraic logics form a Leibniz class.
\end{Theorem}

\begin{proof}
To prove that $\Phi$ is a Leibniz condition, consider $\alpha, \beta \in \class{OR}$ such that $\alpha \leq \beta$. The logic $\vdash_{\class{P}}^{\alpha}$ is protoalgebraic by Proposition \ref{Prop:proto-class1}(ii). This fact and
\[
\vert \LL_{\alpha} \vert = \omega + \vert \alpha \vert \leq \omega + \vert \beta \vert
\]
allow us to apply Proposition \ref{Prop:proto-class1}(i), obtaining ${\vdash_{\class{P}}^{\beta}} \leq {\vdash_{\class{P}}^{\alpha}}$. Hence we conclude that $\Phi$ is a Leibniz condition. Finally, the fact that $\mathrm{Log}(\Phi)$ is the class of protoalgebraic logics is a direct consequence of Proposition \ref{Prop:proto-class1}.
\end{proof}

By Remark \ref{Rem:equivalential} we know that the collection $\class{Equiv}$ of equivalential logics is a Leibniz class. It is therefore sensible to wonder whether we can find an intelligible Leibniz condition $\Phi$ such that $\class{Equiv} = \mathrm{Log}(\Phi)$. This can be done with a simple adaptation of the method employed in the case of protoalgebraic logics.

\begin{law}\label{Def:equiv-basic}
Given an infinite cardinal $\kappa$, let $\mathscr{L}^{\kappa}_{\class{E}}$ be the language consisting of the binary symbols $\{ \multimap_{\alpha} \colon \alpha < \kappa \}$. We set
\[
\Delta_{\kappa}(x, y) \coloneqq \{ x \multimap_{\alpha} y \colon \alpha < \kappa \}.
\]
The \textit{basic equivalential logic of rank $\kappa$} is the logic $\vdash_{\class{E}}^{\kappa}$ formulated on $Fm_{\LL_{\class{E}}^{\kappa}}(\omega)$ determined by the following rules, stipulated for every $\alpha < \kappa$,
\[
\emptyset \, \rhd \, \Delta_{\kappa}(x, x) \qquad x, \Delta_{\kappa}(x, y) \rhd y
\]
\[
\pushQED{\qed}\Delta_{\kappa}(x_{1}, y_{1}) \, \cup \,  \Delta_{\kappa}(x_{2}, y_{2})\rhd \Delta_{\kappa}(x_{1} \multimap_{\alpha} x_{2}, y_{1} \multimap_{\alpha} y_{2}).\qedhere \popQED
\]
\end{law}

The importance of the logic $\vdash_{\class{E}}^{\kappa}$ is justified as follows:

\begin{Proposition}\label{Prop:proto-class1}  Let $\vdash$ be a logic.
\benroman
\item $\vdash$ is equivalential if and only if ${\vdash_{\class{E}}^{\kappa}} \leq {\vdash}$ for every (equiv.\ some) infinite cardinal $\kappa \geq \vert \LL_{\vdash} \vert$.
\item If $\kappa$ is an infinite cardinal and ${\vdash_{\class{E}}^{\kappa}} \leq {\vdash}$, then $\vdash$ is equivalential.
\eroman
\end{Proposition}

\begin{proof}
Condition (ii) follows from the fact that equivalential logics form a Leibniz class by Remark \ref{Rem:equivalential}, while the proof of (i) is analogous to that of Proposition \ref{Prop:proto-class1}(i).
\end{proof}

For every ordinal $\alpha$, let $\vdash_{\class{E}}^{\alpha}$ be the logic $\vdash_{\class{E}}^{\omega + \vert \alpha \vert}$.

\begin{Theorem}\label{Thm:equiv-class}
The sequence $\Phi = \{ \vdash_{\class{E}}^{\alpha}  \colon \alpha \in \class{OR} \}$ is a Leibniz condition and $\mathrm{Log}(\Phi)$ coincides with the class of equivalential logics. In particular, equivalential logics form a Leibniz class.
\end{Theorem}

\begin{proof}
Analogous to the one of Theorem \ref{Thm:proto-class}.
\end{proof}

A logic $\vdash$ is said to be \textit{assertional} \cite{AFRM15,Ra06a} if $F$ is a singleton for every $\langle \A, F \rangle \in \ModS(\vdash)$. 
\begin{Proposition}\label{Prop:asrt-class1}
Let $\vdash$ be a logic.
\benroman
\item  If $x, y, \varphi(x, \vec{z}) \vdash \varphi(y, \vec{z})$ for every formula $\varphi(x, \vec{z}) \in Fm(\vdash)$, then $\vert F \vert \leq 1$ for every $\langle \A, F \rangle \in \ModS(\vdash)$.
\item  $\vdash$ is assertional if and only if it has theorems and $x, y, \varphi(x, \vec{z}) \vdash \varphi(y, \vec{z})$ for every formula $\varphi(x, \vec{z}) \in Fm(\vdash)$.
\eroman
\end{Proposition}

\begin{proof}
See \cite[Thm.\ 10]{AFRM15}. The result is attributed to Suszko in \cite{Cz81}, see also \cite{Rt93}.
\end{proof}
To prove that assertional logics form a \textit{strong} Leibniz class, we need few more concepts.

\begin{law}\label{Def:asrt-basic}
The \textit{basic assertional logic} is the logic $\vdash_{\class{A}}$ formulated in countably many variables and in the language comprising just a unary connective $\top(x)$, axiomatized by the rule $\emptyset \rhd \top(x)$.
\end{law}

A \textit{pointed set} is an algebra $\A = \langle A; \top^{\A} \rangle$ where $\top^{\A}$ is a unary constant map on $A$. In this case, we denote by $\top_\ast^{\A}$ the element of $A$ defined by the map $\top^{\A} \colon A \to A$.

\begin{Proposition}\label{Prop:asrt-class2}
The logic $\vdash_{\class{A}}$ is assertional and
\[
\ModS(\vdash_{\class{A}}) = \{ \langle \A, \{ \top^{\A}_{\ast} \} \rangle \colon \A \text{ is a pointed set}\}.
\]
\end{Proposition}

\begin{proof}
Due to the poor signature of $\vdash_{\class{A}}$ and the fact that $\emptyset \vdash_{\class{A}}\top(x)$, it is easy to see that $\emptyset \vdash_{\class{A}} \varphi$ for every formula $\varphi \in Fm(\vdash_{\class{A}})$ that is not a variable. As a consequence, we obtain $x, y, \varphi(x, \vec{z}) \vdash \varphi(y, \vec{z})$ for every formula $\varphi(x, \vec{z}) \in Fm(\vdash_{\class{A}})$. Moreover, the logic $\vdash_{\class{A}}$ has theorems, e.g.\ $\top(x)$. Hence, we can apply Proposition \ref{Prop:asrt-class1}(ii) obtaining that $\vdash_{\class{A}}$ is assertional.

Now we turn to prove the equality in the statement. First consider a matrix $\langle \A, F \rangle \in \ModS(\vdash_{\class{A}})$. Since $\vdash_{\class{A}}$ is assertional,  $F = \{ a \}$ for some $a \in A$. Together with the fact that $\top(x)$ is a theorem, this yields $\top^{\A}(c) = a$ for all $c \in A$. Hence $\top^{\A}$ is a constant map on $A$ and $a = \top^{\A}_{\ast}$. This establishes the inclusion from left to right.

To prove the other inclusion, let $\A$ be a pointed set. It is clear that the matrix $\langle \A, \{ \top^{\A}_{\ast} \} \rangle$ is a model of the rule $\emptyset \rhd \top(x)$ and, therefore, of $\vdash_{\class{A}}$. It only remains to prove that the congruence $\tarski_{\vdash_{\class{A}}}^{\A} \{ \top_{\ast}^{\A}\}$ is the identity relation. To this end, consider two distinct elements $a, c \in A$. We can assume without loss of generality that $c \ne \top^{\A}_{\ast}$. Then, by the compatibility of the Leibniz congruence, we obtain $\langle a, c \rangle \notin \leibniz^{\A}\{ \top_{\ast}^{\A}, a \}$. Moreover, observe that the matrix $\langle \A, \{ \top_{\ast}^{\A}, a \} \rangle$ is a model of the rule $\emptyset \rhd \top(x)$ and, therefore, of the logic $\vdash_{\class{A}}$, whence $\tarski_{\vdash_{\class{A}}}^{\A} \{ \top_{\ast}^{\A}\} \subseteq \leibniz^{\A}\{ \top_{\ast}^{\A}, a \}$. Hence we conclude that $\langle a, c \rangle \notin \tarski_{\vdash_{\class{A}}}^{\A} \{ \top_{\ast}^{\A}\}$. This shows that $\tarski_{\vdash_{\class{A}}}^{\A} \{ \top_{\ast}^{\A}\}$ is the identity relation and, therefore, $\langle \A, \{ \top_{\ast}^{\A}\} \rangle \in \ModS(\vdash_{\class{A}})$.
\end{proof}

As a consequence, we obtain the following:

\begin{Theorem}\label{Thm:StrongAssertional}
Assertional logics form the strong Leibniz class $\mathrm{Log}(\vdash_{\class{A}})$.
\end{Theorem}

\begin{proof}
It suffices to show that a logic $\vdash$ is assertional if and only if ${\vdash_{\class{A}}} \leq {\vdash}$. To prove the ``if'' part, suppose that ${\vdash_{\class{A}}} \leq {\vdash}$ and let $\btau$ be an interpretation of $\vdash_{\class{A}}$ into $\vdash$. Then consider $\langle \A, F \rangle \in \ModS(\vdash)$. We have that $\langle \A^{\btau}, F\rangle \in \ModS(\vdash_{\class{A}})$. Together with Proposition \ref{Prop:asrt-class2}, this yields that $F$ is a singleton. Hence we conclude that $\vdash$ is assertional.

To prove the ``only if' part'', suppose that $\vdash$ is assertional. From Proposition \ref{Prop:asrt-class1}(ii) we know that $\vdash$ has a theorem $\varphi(x)$. Then let $\btau$ be the translation of $\LL_{\vdash_{\class{A}}}$ into $\LL_{\vdash}$ that sends $\top(x)$ to $\varphi(x)$. We shall see that $\btau$ is an interpretation of $\vdash_{\class{A}}$ into $\vdash$. To this end, consider a matrix $\langle \A, F \rangle \in \ModS(\vdash)$. Since $\vdash$ is assertional, there is $a \in A$ such that $F = \{ a \}$. Together with the fact that $\varphi(x)$ is a theorem of $\vdash$, this implies that $\varphi^{\A}$ is the constant map with value $a$. Hence $\A^{\btau} = \langle A; \varphi^{\A} \rangle$ is essentially a pointed set $\B$, and $\langle \A^{\btau}, F \rangle = \langle \B, \{ \top^{\B}_{\ast}\}\rangle$. By Proposition \ref{Prop:asrt-class2} this guarantees that $\langle \A^{\btau}, F \rangle \in \ModS(\vdash_{\class{A}})$. Hence we conclude that $\btau$ is an interpretation of $\vdash_{\class{A}}$ into $\vdash$.
\end{proof}

A logic $\vdash$ is \textit{truth-equational} \cite{Ra06a} if there is a set $E(x)$ of equations  such that for every $\langle \A, F \rangle \in \ModS(\vdash)$, and $a \in A$,
\[
a \in F \Longleftrightarrow \A \vDash E(a).
\]
Similarly, a logic $\vdash$ is said to be \textit{parametrically truth-equational} \cite{TMo15} if there is a set $E(x, \vec{y})$ of equations  such that for every $\langle \A, F \rangle \in \ModS(\vdash)$ with $F \ne \emptyset$, and $a \in A$,
\[
a \in F \Longleftrightarrow \A \vDash E(a, \vec{c}) \text{ for every }\vec{c} \in A.
\]

\begin{Theorem}\label{Thm:truth-class}
Parametrically truth-equational logics, and truth-equational logics  form  Leibniz classes.
\end{Theorem}

\begin{proof}
We detail only the proof of the fact that parametrically truth-equational logics form a Leibniz class. By Theorem \ref{Thm:Main} it will be enough to show that parametrically truth-equational logics are closed under term-equivalence, compatible expansions, and non-indexed products of sets. The fact that they are closed under term-equivalence, and compatible expansions is clear. Then consider a family $\{  {\vdash_{i}} \colon i \in I \}$ of parametrically truth-equational logics. For every $j \in I$, let $E_{j}(x, \vec{y})$ be the set of equations witnessing the fact that $\vdash_{j}$ is parametrically truth-equational. For every formula $\varphi(x, \vec{y})$ of $\vdash_{j}$, we denote by $\hat{\varphi}$ the sequence $\langle \varphi_{i} \colon i \in I \rangle$ in which $\varphi_{i} = x$ for every $i \in I \smallsetminus \{ j \}$, and $\varphi_{j} = \varphi$. Observe that $\hat{\varphi}(x, \vec{y})$ is a basic operation of $\bigotimes_{i \in I}{\vdash_{i}}$. Bearing this in mind, we define
\[
E(x, \vec{y}) \coloneqq \{ \hat{\varphi} \thickapprox \hat{\psi} \colon \varphi \thickapprox \psi \in \bigcup_{i \in I}E_{i}(x, \vec{y}) \}.
\]
Then consider a matrix $\langle \A, F \rangle \in \ModS(\bigotimes_{i \in I}{\vdash_{i}})$ such that $F \ne \emptyset$. By \cite[Prop.\ 4.5]{JaMor19-1} we have  $\langle \A, F \rangle \leq_{sd} \bigotimes_{i \in I}\langle \A_{i}, F_{i}\rangle$ for some $\langle \A_{i}, F_{i}\rangle \in \ModS(\vdash_{i})$ with $F_{i} \ne \emptyset$. For every $a \in A$,
\begin{align*}
a \in F &\Longleftrightarrow a(i) \in F_{i} \text{ for every }i \in I\\
&\Longleftrightarrow \A_{i} \vDash E_{i}(a(i), \vec{c}) \text{ for every }\vec{c} \in A_{i} \text{ and }i \in I\\
&\Longleftrightarrow \A \vDash E(a, \vec{c}) \text{ for every }\vec{c} \in A.
\end{align*}
The above equivalences are justified as follows: the first is straightforward, the second follows from the fact the set $E_{i}(x, \vec{y})$ witnesses that $\vdash_{i}$ is parametrically truth-equational, and the third is a consequence of the fact that the projection $\pi_{i} \colon A \to A_{i}$ is surjective for all $i \in I$. From the display above we obtain that $\bigotimes_{i \in I}{\vdash_{i}}$ is parametrically truth-equational, as desired.
\end{proof}

\begin{problem}
Is it possible to find a transparent Leibniz condition $\Phi$ such that $\mathrm{Log}(\Phi)$ is the class of (parametrically) truth-equational logics?
\end{problem}

A logic $\vdash$ is said to be \textit{order algebraizable} \cite{JRa13} if there is a set $\Delta(x, y)$ of formulas  and a set $E(x)$ of inequalities such that for every $\langle \A, F \rangle \in \ModS(\vdash)$ the relation $\preccurlyeq_{F}^{\A}$ on $A$ defined as follows is a partial order: for every $a, c \in A$,
\[
a \preccurlyeq_{F}^{\A}c \Longleftrightarrow \Delta^{\A}(a, c) \subseteq F,
\]
and for every $a \in A$,
\[
a \in F \Longleftrightarrow \langle \A, \preccurlyeq_{F}^{\A} \rangle \vDash E(e).
\]

\begin{Theorem}\label{Thm:order-class}
Order algebraizable logics, and logics with theorems form Leibniz classes.
\end{Theorem}

\begin{proof}
The result can be established by checking condition (ii) of Theorem \ref{Thm:Main} with ideas similar to the ones in the proof of Theorem \ref{Thm:truth-class}.
\end{proof}

We expect that a Leibniz condition defining the class of order algebraizable logics could be extracted from \cite[Thm.\ 7.1(ii)]{JRa13}.

\begin{Remark}\label{Rem:algebraizable-logics}
Recall from Proposition \ref{Prop:meets-and-joins} that every class of logics that can be written as the intersection of two Leibniz classes is still a Leibniz class. This observation can be exploited to show that some well-known collections of logics are Leibniz classes. For instance, a logic $\vdash$ is said to be \textit{algebraizable} \cite{BP89} (resp.\ \textit{weakly algebraizable} \cite{CzJa00}) if is equivalential (resp.\ protoalgebraic) and truth-equational. From Theorems \ref{Thm:proto-class}, \ref{Thm:truth-class}, and \ref{Thm:equiv-class} it follows that (weakly) algebraizable logics form a Leibniz class.
\qed
\end{Remark}

We conclude this section by providing some examples of collections of logics that are not Leibniz classes. By Remark \ref{Rem:extensions} we know that every class of logics that is not closed under the formation of extensions will serve this purpose. Among these we count, for instance, the class of logics with the Craig deductive interpolation property \cite{CzP99}, and the class of logics with the infinite (resp.\ finite) Beth definability property \cite{BMR16,BlHoo06}. 

To describe another interesting class of logics that is not a Leibniz class, recall that a logic $\vdash$ has an \textit{algebraic semantics} \cite{BP89,BlRe03} if there are a set of equations $E(x)$ and a class of algebras $\class{K}$ such that for every $\Gamma \cup \{ \varphi \} \subseteq Fm(\vdash)$,
\[
\Gamma \vdash \varphi \Longleftrightarrow \bigcup \{ E(\gamma) \colon \gamma \in \Gamma \} \models_{\class{K}} E(\varphi)
\]
where $\models_{\class{K}}$ is the equational consequence relative to $\class{K}$. Surprisingly enough, every logic is term-equivalent to one with an algebraic semantics, as shown essentially in \cite[Thm.\ 3.1]{BlRe03}. Together with the fact that some logics lack an algebraic semantics \cite{BK83}, this yields that the class of logics with an algebraic semantics is not closed under term-equivalence, whence it is not a Leibniz class by Theorem \ref{Thm:Main}(ii). 

The reader familiar with abstract algebraic logic may be interested to know that also (fully) selfextensional, and (fully) Fregean logics \cite{FJa09} do not form Leibniz classes, since these collections are not closed under compatible expansions.

\begin{problem}
The majority of well-known Leibniz classes can be characterized in terms of the behaviour of the so-called \textit{Leibniz operator}, i.e.\ the map $\leibniz^{\A} \colon \mathcal{P}(A) \to \Con\A$, defined on every algebra $\A$, that associates a subset $F \subseteq A$ with the Leibniz congruence $\leibniz^{\A}F$. Is it possible to establish a precise relation between Leibniz classes and the behaviour of the Leibniz operator?
\end{problem}

Later on we make use of the following well-known observation \cite{AAL-AIT-f}.

\begin{Proposition}\label{Prop:still-have-thms}
The classes of protoalgebraic, equivalential, assertional, order algebraizable, and truth-equational logics comprise only logics with theorems.
\end{Proposition}

\section{Irreducibility, primality and their boundaries}

Part of the speculative power of the identification of the Leibniz hierarchy with the lattice of all Leibniz classes comes from the fact that it allows to apply order-theoretic methods and intuitions to the study of the first. To explain how, recall that an element $a$ of a lattice $\langle A; \land, \lor \rangle$ is said to be \textit{meet-irreducible} if for every pair $b, c \in A$,
\[
\text{if }a = b \land c\text{, then either }a = b\text{ or }a = c.
\]
Similarly, $a$ is said to be \textit{meet-prime} if for every pair $b, c \in A$,
\[
\text{if }b \land c \leq a \text{, then either }b \leq a\text{ or }c \leq a.
\]
Accordingly, $a$ is said to be \textit{meet-reducible} when it is not meet-irreducible. It is clear that every meet-prime element of $\langle A; \land, \lor \rangle$ is meet-irreducible, while the converse is not true in general. 

Since the Leibniz hierarchy is a lattice (Proposition \ref{Prop:meets-and-joins}), it makes sense to ask whether a Leibniz class is meet-irreducible or meet-prime. An affirmative answer to this question can then be regarded as a certificate that the Leibniz class under consideration captures a primitive or fundamental concept.

\begin{Remark}
Since infima in the Leibniz hierarchy are intersections, a Leibniz class $\oper{K}$ is meet-prime if there is no pair of logics ${\vdash_{1}}, {\vdash_{2}} \notin \oper{K}$ such that ${\vdash} \in \oper{K}$, for every logic $\vdash$ such that ${\vdash_{1}}, {\vdash_{2}} \leq {\vdash}$. 
\qed
\end{Remark}

\begin{exa}
The Leibniz class of algebraizable (resp.\ weakly algebraizable) logics is meet-reducible, since it can be obtained as the intersection of the strictly larger Leibniz classes of equivalential (resp.\ protoalgebraic) logics and truth-equational logics (Remark \ref{Rem:algebraizable-logics}). Moreover, there are meet-irreducible Leibniz classes that are not meet-prime, e.g.\ the class of logics $\vdash$ for which there is no three-element algebra $\A$ and $a \in A$ such that $\langle \A, \{ a \} \rangle \in \ModS(\vdash)$. Even if we do not pursue the details here, the proof of this fact is an adaptation of an argument in \cite[pag.\ 54]{GaTa84}.
\qed
\end{exa}

The next results put some boundaries to the expectation that well-known Leibniz classes should be meet-irreducible or meet-prime.

\begin{Proposition}\label{Prop:dec-among-thm}
All Leibniz classes properly included into the class of logics with theorems are meet-reducible.
\end{Proposition}

\begin{proof}[Proof sketch.]
Given a logic $\vdash$ with theorems, we denote by $\vdash_{\emptyset}$ the logic on $Fm(\vdash)$ defined for every $\Gamma \cup \{ \varphi \} \subseteq Fm(\vdash)$ as 
\[
\Gamma \vdash_{\emptyset} \varphi \Longleftrightarrow \Gamma \ne \emptyset \text{ and }\Gamma \vdash \varphi.
\]
It is clear that $\vdash_{\emptyset}$ lacks theorems. It is not hard to see that
\begin{equation}\label{Eq:empty-logic}
\ModS(\vdash_{\emptyset}) = \ModS(\vdash) \cup \{ \langle \A, \emptyset \rangle \colon \langle \A, F \rangle \in \ModS(\vdash) \text{ for some }F \subseteq A \}.
\end{equation}

Then consider a Leibniz class $\oper{K}$ properly included into the Leibniz class $\class{Thrms}$ of logics with theorems. We set
\[
\oper{K}_{\emptyset} \coloneqq \oper{K} \cup \{ {\vdash_{\emptyset}} \colon {\vdash} \in \oper{K} \}.
\]
It is clear that $\oper{K} \subsetneq \oper{K}_{\emptyset}$, since $\oper{K}_{\emptyset}$ contains logics without theorems. Moreover, $\oper{K} = \oper{K}_{\emptyset} \cap \class{Thrms}$. Therefore, to conclude that $\oper{K}$ is meet-reducible, it will be enough to show that $\oper{K}_{\emptyset}$ is a Leibniz class.

To prove this, consider a Leibniz condition $\Phi = \{ {\vdash^{\alpha}} \colon \alpha \in \class{OR} \}$ such that $\oper{K} = \mathrm{Log}(\Phi)$. We define $\Phi_{\emptyset} \coloneqq \{ {\vdash^{\alpha}_{\emptyset}} \colon \alpha \in \class{OR} \}$. From the fact that $\Phi$ is a Leibniz condition and (\ref{Eq:empty-logic}) it follows that $\Phi_{\emptyset}$ is also Leibniz condition. We shall prove that $\oper{K}_{\emptyset} = \mathrm{Log}(\Phi_{\emptyset})$. To this end, consider a logic ${\vdash} \in \oper{K}$. Then there is $\alpha \in \class{OR}$ such that ${\vdash^{\alpha}} \leq {\vdash}$. Together with (\ref{Eq:empty-logic}), this implies that ${\vdash_{\emptyset}^{\alpha}} \leq {\vdash}, {\vdash_{\emptyset}}$. Hence we obtain $\oper{K}_{\emptyset} \subseteq \mathrm{Log}(\Phi_{\emptyset})$. To prove the other inclusion, consider $\alpha \in \class{OR}$ and a logic $\vdash$ with an interpretation $\btau$ of $\vdash_{\emptyset}^{\alpha}$ into $\vdash$. If $\vdash$ has theorems, then we can use (\ref{Eq:empty-logic}) to conclude that $\btau$ is also an interpretation of $\vdash^{\alpha}$ into $\vdash$, whence ${\vdash} \in \mathrm{Log}(\Phi) = \oper{K} \subseteq \oper{K}_{\emptyset}$. Then we consider the case where $\vdash$ lacks theorems. Since ${\vdash^{\alpha}} \in \oper{K}$, the logic $\vdash^{\alpha}$ has a theorem $\varphi(x)$. Let $y$ be a variable different from $x$ and observe that $y \vdash_{\emptyset}^{\alpha} \varphi(x)$. By \cite[Prop.\ 3.3]{JaMor19-1} this yields $y \vdash \btau(\varphi(x))$. Then consider the logic $\vdash^{+}$ on $Fm(\vdash)$ induced by the class of matrices
\begin{equation}\label{Eq:empty-logic2}
\{ \langle \A, F \rangle \in \ModS(\vdash) \colon F \ne \emptyset \}.
\end{equation}
Since $y \vdash \btau(\varphi(x))$, it is clear that $\emptyset \vdash^{+} \btau(\varphi)$.  In particular, this guarantees that $\ModS(\vdash^{+})$ is the class of matrices in (\ref{Eq:empty-logic2}). This fact, together with (\ref{Eq:empty-logic}), yields that $\btau$ is an interpretation of $\vdash^{\alpha}$ into $\vdash^{+}$. As a consequence, we obtain ${\vdash^{+}} \in \oper{K}$. Since $\vdash$ coincides with $\vdash^{+}_{\emptyset}$, we conclude that ${\vdash} \in \oper{K}_{\emptyset}$. This establishes that $\oper{K}_{\emptyset} = \mathrm{Log}(\Phi_{\emptyset})$ and, therefore, that $\oper{K}$ is a Leibniz class.
\end{proof}

\begin{Remark}
Proposition \ref{Prop:dec-among-thm} indicates that Leibniz classes comprising only logics with theorems (but not all of them) cannot be meet-irreducible in an absolute sense. For this reason, we say that a Leibniz class $\oper{K}$ is meet-reducible \textit{among logics with theorems} if for every pair of Leibniz classes $\oper{K}_{1}$ and $\oper{K}_{2}$ comprising logics with theorems only, if $\oper{K} = \oper{K}_{1} \cap \oper{K}_{2}$, then either $\oper{K} = \oper{K}_{1}$ or $\oper{K} = \oper{K}_{2}$. A similar definition applies to the case of meet-prime Leibniz classes.
\qed
\end{Remark}

\begin{Proposition}
All Leibniz classes properly included into the class of assertional logics are not meet-prime among Leibniz classes with theorems.
\end{Proposition}

\begin{proof}
Let  $\oper{F}$  be the class of logics $\vdash$ with theorems such that if $\langle \A, F \rangle \in \ModS(\vdash)$, then either the matrix $\langle \A, F \rangle$ is trivial or $\vert F \vert \geq 2$. Using the characterization of Leibniz classes given in Theorem \ref{Thm:Main}(ii) it is not hard to see that $\oper{F}$ is indeed a Leibniz class. 

Now, let $\class{Asrt}$ be the Leibniz class of assertional logics, and consider an arbitrary Leibniz class  $\oper{K}$ properly included into $\class{Asrt}$. It is clear that both $\class{Asrt}$ and $\oper{F}$ are not included into $\oper{K}$. Bearing this in mind, it only remains to show that $\class{Asrt} \cap \oper{F} \subseteq \oper{K}$. To this end, consider a logic ${\vdash} \in \class{Asrt} \cap \oper{F}$. Observe that every matrix $\langle \A, F \rangle \in \ModS(\vdash)$ is such that $F$ is a singleton (as ${\vdash} \in \class{Asrt}$) and, therefore, is trivial (as ${\vdash} \in \oper{F}$). Thus $\ModS(\vdash)$ is the class of trivial matrices in the language of $\LL_{\vdash}$. As a consequence, the logic $\vdash$ is inconsistent \cite[Lem.\ 7.1]{JaMor19-1}. In particular, this guarantees that $\llbracket \vdash \rrbracket$ is the maximum of $\class{Log}$  \cite[Thm.\ 7.3]{JaMor19-1}. By Theorem \ref{Thm:Main}(iii) the collection $\oper{K}^{\dagger}$ is a non-empty filter of $\class{Log}$ which implies  $\llbracket \vdash \rrbracket \in \oper{K}^{\dagger}$. Together with $\oper{K} = \{ \vdash' \colon \llbracket \vdash' \rrbracket \in \oper{K}^{\dagger}\}$, this yields ${\vdash} \in \oper{K}$, as desired.
\end{proof}

\section{Meet-prime Leibniz classes}

In this section we show that the Leibniz class  of all logics with theorems, and the Leibniz class of truth-minimal logics that we introduce below are meet-prime in the absolute sense. On the other hand, it is proved that the Leibniz classes of assertional and truth-equational logics are meet-prime among the logics with theorems. As both assertional and truth-equational logics have theorems, Proposition \ref{Prop:dec-among-thm} guarantees that the restriction to logics with theorems cannot be dropped here. For the present purpose, it is convenient to start the discussion from the new class of truth-minimal logics.

\subsection{Truth-minimal logics}

\begin{law}
A logic $\vdash$ is \textit{truth-minimal} if for every $\langle\A, F\rangle \in \ModS(\vdash)$,
\[
\text{if }\langle \A, G\rangle \in \ModS(\vdash) \text{ and }G \subseteq F\text{, then either }G = F\text{ or }G = \emptyset.
\]
\end{law}

The next result is instrumental to construct examples of truth-minimal logics.

\begin{Proposition}\label{Prop:min-trick3}
Let $\kappa$ be an infinite cardinal and $\class{M}$ a class of matrices such that $\vert F \vert \leq 1$ for every $\langle \A, F \rangle \in \class{M}$. The logic $\vdash$  induced by $\class{M}$ formulated in $\kappa$ variables is truth-minimal.
\end{Proposition}

\begin{proof}
The assumption that $\vert F \vert \leq 1$ for every $\langle \A, F \rangle \in \class{M}$ implies that $x, y, \varphi(x, \vec{z}) \vdash \varphi(y, \vec{z})$ for every formula $\varphi(x, \vec{z}) \in Fm(\vdash)$. By Proposition \ref{Prop:asrt-class1}(i) we conclude that $\vdash$ is truth-minimal.
\end{proof}

Our first aim is to show that truth-minimal logics constitute a Leibniz class. To this end, we make use of the following observation:

\begin{Lemma}\label{Lem:min-trick1}
Let $\{ \vdash_{i} \colon i \in I \}$ be a family of logics, and $\A \leq_{\textup{sd}} \bigotimes_{i \in I}\A_{i}$ where $\A_{i}$ is an $\LL_{\vdash_{i}}$-algebra for every $i \in I$. If $G \ne \emptyset$ is a deductive filter of $\bigotimes_{i \in I}{\vdash_{i}}$ on $\A$, then for every $j \in I$ there is a deductive filter $G_{j}$ of $\vdash_{j}$ on $\A_{j}$ such that $G = A \cap \prod_{i \in I}G_{i}$. 
\end{Lemma}

\begin{proof}
For every $j \in I$, we denote the natural projection map by $\pi_{j} \colon \A \to \A_{j}$, and set $G_{j} \coloneqq \pi_{j}[G]$. 

Then we turn to prove that $G = A \cap \prod_{i \in I}G_{i}$. The inclusion from left to right is clear. To prove the other one, consider an element $a \in A \cap \prod_{i \in I}G_{i}$. For every $j \in I$ there is $c_{j} \in G$ such that $c_{j}(j) = a(j)$. Then consider the basic operation $x \multimap_{j} y$ of $\bigotimes_{i \in I}{\vdash_{i}}$, whose $j$-th component is the projection on the first coordinate, and whose $i$-th component is  the projection on the second coordinate for every $i \in I \smallsetminus \{ j \}$. Since $c_{j}(j) = a(j)$,
\begin{equation}\label{Eq:min-1}
a \multimap_{j}^{\A} c_{j} = c_{j} \in G.
\end{equation}

Now, consider some distinct variables $\{ y_{j} \colon j \in I \} \cup \{ x \} \subseteq Fm(\bigotimes_{i \in I}{\vdash_{i}})$. Bearing in mind that $\bigotimes_{i \in I}{\vdash_{i}}$ is the logic induced by the class $\bigotimes_{i \in I}\ModS(\vdash_{i})$, it is easy to see that the following rule is valid in $\bigotimes_{i \in I}{\vdash_{i}}$:
\[
\{ x \multimap_{j} y_{j} \colon j \in I \} \rhd x.
\]
Together with (\ref{Eq:min-1}) and the fact that $G$ is a deductive filter of $\bigotimes_{i \in I}{\vdash_{i}}$ on $\A$, this implies  $a \in G$. Hence we conclude
\begin{equation}\label{Eq:min-3}
G = A \cap \prod_{i \in I}G_{i}.
\end{equation}

Consider an index $j \in I$. To conclude the proof, it suffices to show that $G_{j}$ is a filter of $\vdash_{j}$ on $\A_{j}$. To prove this, suppose that $\Gamma \vdash_{j} \varphi$ and consider a homomorphism $h \colon \Fm(\vdash_{j}) \to \A_{j}$ such that $h[\Gamma] \subseteq G_{j}$. For every $n$-ary term $\psi(z_{1}, \dots, z_{n}) \in Fm(\vdash_{j})$, let $\hat{\psi}(z_{1}, \dots, z_{n})$ be an arbitrary basic $n$-ary operation of $\bigotimes_{i \in I}{\vdash_{i}}$ whose $j$-th component is $\psi$. Moreover, consider a variable  $y \in Fm(\bigotimes_{i \in I}{\vdash_{i}})$ not occurring in $\Gamma$.\footnote{In case there is no such a variable $y$, we replace in the proof $\Gamma \cup \{ \varphi \}$ by $\sigma[\Gamma \cup \varphi]$, where $\sigma$ is a substitution on $Fm(\vdash_{j})$ that maps variables to variables, and that is injective but not surjective.} Again bearing in mind that $\bigotimes_{i \in I}{\vdash_{i}}$ is the logic induced by the class $\bigotimes_{i \in I}\ModS(\vdash_{i})$, it is clear that the following rule is valid in $\bigotimes_{i \in I}{\vdash_{i}}$:
\begin{equation}\label{Eq:min-2}
\{ \hat{\gamma} \multimap_{j} y  \colon \gamma \in \Gamma \} \rhd \hat{\varphi} \multimap_{j} y.
\end{equation}

Now, choose an element $a \in G$, and let $f \colon \Fm(\bigotimes_{i \in I}{\vdash_{i}}) \to \A$ be a homomorphism such that $f(y) = a$ and $f(x)(j) = h(x)$ for every variable $x \in Fm(\vdash_{j})$. Observe that there exists such an $f$, since the projection $\pi_{j} \colon A \to A_{j}$ is surjective. For every $i \in I \smallsetminus \{ j \}$ and $\gamma \in \Gamma$, we have
\begin{align*}
f(\hat{\gamma} \multimap_{j} y)(i) &= a(i) \in \pi_{i}[G] = G_{i}\\
f(\hat{\gamma} \multimap_{j} y)(j) &= h(\gamma)\in G_{j}.
\end{align*}
Together with (\ref{Eq:min-3}), this implies
\[
f(\hat{\gamma} \multimap_{j} y) \in A \cap \prod_{i \in I}G_{i} = G \text{, for every }\gamma \in \Gamma.
\]
The above display, together with (\ref{Eq:min-2}) and the fact that $G$ is a deductive filter of $\bigotimes_{i \in I}{\vdash_{i}}$ on $\A$, yields that $f(\hat{\varphi} \multimap_{j}y) \in G$. As a consequence, we obtain
\[
h(\varphi) = f(\hat{\varphi} \multimap_{j}y)(j) \in \pi_{j}[G] = G_{j}.
\]
Hence we conclude that $G_{j}$ is a deductive filter of $\vdash_{j}$ on $\A_{j}$.
\end{proof}

As a consequence, we obtain the following:

\begin{Theorem}\label{Thm:truth-minimal-class}
The class of truth-minimal logics is a Leibniz class.
\end{Theorem}

\begin{proof}
It is straightforward that the class of truth-minimal logics is closed under term-equivalence and compatible expansions. In the light of Theorem \ref{Thm:Main}, it only remains to prove that this class is closed under the formation of non-indexed products of sets. 

To prove this, consider a family $\{ {\vdash_{i}} \colon i \in I \}$ of truth-minimal logics. Moreover, consider $\langle \A, F \rangle, \langle \A, G \rangle \in \ModS(\bigotimes_{i \in I} {\vdash_{i}})$ with $\emptyset \ne G \subseteq F$. We need to show that $F = G$. By \cite[Prop.\ 4.5]{JaMor19-1} there is a family of matrices $\{ \langle \A_{i}, F_{i}\rangle \in \ModS(\vdash_{i}) \colon i \in I \}$ such that $\langle \A, F \rangle \leq_{\textup{sd}} \bigotimes_{i \in I}\langle \A_{i}, F_{i}\rangle$. By Lemma \ref{Lem:min-trick1}, for every $j \in I$ there is a deductive filter $G_{j}$ of $\vdash_{j}$ on $\A_{j}$ such that
\begin{equation}\label{Eq:minimality-z1}
\emptyset \ne G_{j} \subseteq F_{j}
\end{equation}
and $G = A \cap \prod_{i \in I}G_{i}$. Now, from the fact that $G_{j}$ is a deductive filter of $\vdash_{j}$ and $G_{j} \subseteq F_{j}$ it follows that $\tarski_{\vdash_{j}}^{\A} G_{j} \subseteq \tarski_{\vdash_{j}}^{\A} F_{j}$. Therefore, bearing in mind that $\tarski_{\vdash_{j}}^{\A} F_{j}$ is the identity relation, we conclude that the same holds for $\tarski_{\vdash_{j}}^{\A} G_{j}$. In particular, this yields
\begin{equation}\label{Eq:minimality-z2}
\langle \A_{j}, G_{j}\rangle \in \ModS(\vdash_{j}).
\end{equation}

Finally, for every $j \in J$ we can apply the fact that $\vdash_{j}$ is truth-minimal to (\ref{Eq:minimality-z1}, \ref{Eq:minimality-z2}), obtaining that $F_{j} = G_{j}$. As a consequence, we get
\[
G = A \cap \prod_{i \in I}G_{i} = A \cap \prod_{i \in I}F_{j} = F.
\]
Hence we conclude that $\bigotimes_{i \in I}{\vdash_{i}}$ is truth-minimal.
\end{proof}

To prove that the Leibniz class of truth-minimal logics is meet-prime, we rely on the following technical observations:

\begin{Lemma}\label{Lem:truth-prime0}
Let $\vdash$ be a logic. If $\langle \A, A \rangle \in \ModS(\vdash)$, then $\A$ is the trivial algebra.
\end{Lemma}

\begin{proof}
From \cite[Prop.\ 2.2(ii)]{JaMor19-1} we obtain $\tarski_{\vdash}^{\A}A = A \times A$. On the other hand, since $\langle \A, A \rangle \in \ModS(\vdash)$, the congruence $\tarski_{\vdash}^{\A}A$ is the identity relation. But the fact that $A \times A$ is the identity relation implies that $A$ is a singleton, as desired.
\end{proof}

\begin{Lemma}\label{Lem:truth-prime1}
Let $\vdash$ be a logic, and $\langle \A, F \rangle, \langle \A, G \rangle \in \ModS(\vdash)$ such that $\emptyset \ne G \subsetneq F$. There are $\langle \B, F^{\dagger} \rangle, \langle \B, G^{\dagger} \rangle \in \ModS(\vdash)$ such that $\emptyset \ne G^{\dagger} \subsetneq F^{\dagger} \subsetneq B$ and
\[
\omega \leq \vert B \vert = \vert G^{\dagger} \vert = \vert F^{\dagger} \smallsetminus G^{\dagger} \vert = \vert B \smallsetminus F^{\dagger} \vert.
\]
\end{Lemma}

\noindent \textit{Proof.} Choose an infinite cardinal $\kappa \geq \vert A \vert$ and define $\B \coloneqq \A^{\kappa} \times \A^{\kappa}$. Then set
\[
F^{\dagger} \coloneqq F^{\kappa} \times F^{\kappa} \text{ and }G^{\dagger} \coloneqq F^{\kappa} \times G^{\kappa}.
\]
Since $\ModS(\vdash)$ is closed under direct products \cite[Lem.\ 2.3]{JaMor19-1}, we obtain
\[
\langle \B, F^{\dagger}\rangle, \langle \B, G^{\dagger}\rangle \in \ModS(\vdash).
\]
Moreover, from $\emptyset \ne G \subsetneq F$ it follows $\emptyset \ne G^{\dagger} \subsetneq F^{\dagger}$.

Observe that since $\emptyset \ne G \subsetneq F$, the set $F$ contains at least two distinct elements. This fact and $\emptyset \ne G \subsetneq F$ guarantee that
\begin{align}
2^{\kappa} &\leq \vert F^{\kappa} \times G^{\kappa} \vert = \vert G^{\dagger} \vert \label{Eq:truth-prime2} \\
2^{\kappa} &\leq \vert F^{\kappa} \vert \leq \vert F^{\kappa} \times (F^{\kappa} \smallsetminus G^{\kappa}) \vert = \vert F^{\dagger}  \smallsetminus G^{\dagger} \vert\label{Eq:truth-prime3}
\end{align}
Since $F$ contains at least two elements and $\langle \A, F \rangle \in \ModS(\vdash)$, we can apply Lemma \ref{Lem:truth-prime0} obtaining $F \subsetneq A$. Bearing this in mind, we get $F^{\dagger} \subsetneq B$ and
\begin{equation}\label{Eq:truth-prime4}
2^{\kappa} \leq \vert A^{\kappa} \vert \leq \vert A^{\kappa} \times (A^{\kappa} \smallsetminus F^{\kappa}) \vert \leq \vert B \smallsetminus F^{\dagger} \vert.
\end{equation}

Finally, from $\kappa \geq \omega + \vert A \vert$ it follows  $\vert B \vert = \vert A^{\kappa} \times A^{\kappa} \vert = 2^{\kappa}$. But, together with (\ref{Eq:truth-prime2}, \ref{Eq:truth-prime3}, \ref{Eq:truth-prime4}), this yields
\[
\pushQED{\qed} \omega \leq \vert B \vert = \vert G^{\dagger} \vert = \vert F^{\dagger} \smallsetminus G^{\dagger} \vert = \vert B \smallsetminus F^{\dagger} \vert.\qedhere \popQED
\]

\begin{Lemma}\label{Lem:truth-prime2}
Let $\vdash_{1}$ and $\vdash_{2}$ be logics, and for every $i = 1, 2$ let $\langle \A_{i}, F_{i} \rangle, \langle \A_{i}, G_{i} \rangle \in \ModS(\vdash_{i})$ such that $\emptyset \ne G_{i} \subsetneq F_{i}$. There are sets $\emptyset \ne G \subsetneq F \subsetneq B$, an $\LL_{\vdash_{1}}$-algebra $\B_{1}$, and an $\LL_{\vdash_{2}}$-algebra $\B_{2}$ such that $B_{1} = B_{2} = B$, and 
\begin{align*}
\langle \B_{1}, F \rangle, \langle \B_{1}, G \rangle &\in \ModS(\vdash_{1})\\
\langle \B_{2}, F \rangle, \langle \B_{2}, G \rangle &\in \ModS(\vdash_{2}).
\end{align*}
\end{Lemma}

\begin{proof} 
By Lemma \ref{Lem:truth-prime1} for every $i = 1, 2$ there are $\langle \C_{i}, F^{\dagger}_{i} \rangle, \langle \C_{i}, G^{\dagger}_{i} \rangle \in \ModS(\vdash_{i})$ such that $\emptyset \ne G^{\dagger}_{i} \subsetneq F^{\dagger}_{i} \subsetneq C_{i}$ and $\omega \leq \vert C_{i} \vert = \vert G^{\dagger}_{i} \vert = \vert F^{\dagger}_{i} \smallsetminus G^{\dagger}_{i} \vert = \vert C_{i} \smallsetminus F^{\dagger}_{i} \vert$.

 Let then $\kappa \coloneqq \max \{ \vert C_{1} \vert, \vert C_{2} \vert \}$. Consider also $i = 1, 2$. Since $\ModS(\vdash_{i})$ is closed under direct powers \cite[Lem.\ 2.3]{JaMor19-1}, we obtain
\begin{equation}\label{Eq:truth-prime5}
\langle \C_{i}^{\kappa}, (F^{\dagger}_{i})^{\kappa} \rangle, \langle \C_{i}^{\kappa}, (G^{\dagger}_{i})^{\kappa} \rangle \in \ModS(\vdash_{i}).
\end{equation}
Moreover, since $\emptyset\ne G_{i}^{\dagger} \subsetneq F^{\dagger}_{i} \subsetneq C_{i}$, we have
\[
\emptyset\ne G_{i}^{\dagger\kappa} \subsetneq F^{\dagger\kappa}_{i} \subsetneq C_{i}^{\kappa}.
\]
Finally, since $\omega \leq \vert C_{i} \vert = \vert G^{\dagger}_{i} \vert = \vert F^{\dagger}_{i} \smallsetminus G^{\dagger}_{i} \vert = \vert C_{i} \smallsetminus F^{\dagger}_{i} \vert \leq \kappa$, we have
\begin{equation}\label{Eq:truth-prime6}
\kappa = \vert C_{i}^{\kappa} \vert = \vert G^{\dagger\kappa}_{i} \vert = \vert F^{\dagger\kappa}_{i} \smallsetminus G^{\dagger\kappa}_{i} \vert = \vert C_{i}^{\kappa} \smallsetminus F^{\dagger\kappa}_{i} \vert.
\end{equation}

Then consider some sets $\emptyset \ne G \subsetneq F \subsetneq B$ such that
\[
\kappa = \vert B \vert = \vert G \vert = \vert F \smallsetminus G \vert = \vert B \smallsetminus F \vert.
\]
From (\ref{Eq:truth-prime5}) and (\ref{Eq:truth-prime6}) it easily follows that for each $i = 1, 2$ there is an $\LL_{\vdash_{i}}$-algebra $\B_{i}$ with universe $B$ such that $\langle \B_{i}, F \rangle \cong \langle \C_{i}^{\kappa}, (F^{\dagger}_{i})^{\kappa} \rangle$, $\langle \B_{i}, G\rangle 
\cong \langle \C_{i}^{\kappa}, (G^{\dagger}_{i})^{\kappa} \rangle$ and, therefore, $\langle \B_{i}, F \rangle, \langle \B_{i}, G\rangle \in \ModS(\vdash_{i})$.
\end{proof}

We are now ready to prove the main result of this part:

\begin{Theorem}\label{Thm:truth-min-prime} \
The Leibniz class of truth-minimal logics is meet-prime.
\end{Theorem}

\begin{proof}
Consider two logics $\vdash_{1}$ and $\vdash_{2}$ that are not truth-minimal. It will be enough to construct a logic $\vdash$ that is not truth-minimal and in which $\vdash_{1}$ and $\vdash_{2}$ are interpretable.

To this end, observe that for every $i = 1, 2$ there are $\langle \A_{i}, F_{i}\rangle, \langle \A_{i}, G_{i}\rangle \in \ModS(\vdash_{i})$ such that $\emptyset \ne G_{i} \subsetneq F_{i}$. From Lemma \ref{Lem:truth-prime2} we obtain 
\[
\langle \B_{1}, F \rangle, \langle \B_{1}, G \rangle \in \ModS(\vdash_{1}) \text{ and }\langle \B_{2}, F \rangle, \langle \B_{2}, G \rangle \in \ModS(\vdash_{2})
\]
for some sets $\emptyset \ne G \subsetneq F \subsetneq B$, an $\LL_{\vdash_{1}}$-algebra $\B_{1}$, and an $\LL_{\vdash_{2}}$-algebra $\B_{2}$ such that $B_{1} = B_{2} = B$. 

Let $\B$ be the common expansion of $\B_{1}$ and $\B_{2}$ with all finitary operations on $B$. Moreover, let $\vdash$ be the logic, formulated in $\kappa \coloneqq \max\{ \vert Fm(\vdash_{1}) \vert, \vert Fm(\vdash_{2}) \vert \}$ variables induced by the set of matrices 
\[
\class{K} \coloneqq \{ \langle \B, F \rangle, \langle \B, G \rangle\}. 
\]

We claim that $\vdash$ is equivalential. To prove this, consider elements $0 \in B \smallsetminus F$ and $1 \in G$ (this is possible, since $\emptyset \ne G$ and  $F\subsetneq B$). Then let $x \multimap y$ be the binary basic operation of $\B$ defined as follows for every $a, c \in B$,
\begin{displaymath}
a \multimap^{\B} c \coloneqq \left\{\begin{array}{@{\,}ll}
1 & \text{if $a=c$}\\
0 & \text{otherwise.}\\
\end{array} \right.
\end{displaymath}
Bearing in mind that $0 \notin F \cup G$ and $1 \in F \cap G$, it is not hard to see that the set $\Delta(x, y) \coloneqq \{ x \multimap y \}$ satisfies the conditions in the right hand side of \cite[Thm.\ 2.7]{JaMor19-1}. As a consequence we obtain that $\vdash$ is equivalential, establishing the claim.

Since $\B$ is endowed with all unary constant maps, we have $\SSS(\class{K}) = \class{K}$. Moreover, for each pair of different $a, c \in B$, we have
\[
a \multimap^{\B} c = 0 \notin G \cup F \text{ and }a \multimap^{\B} a = 1 \in G \cap F.
\]
By \cite[Prop.\ 2.2(i)]{JaMor19-1} this implies that the matrices in $\class{K}$ are reduced. This fact, together with the claim and $\SSS(\class{K}) = \class{K}$, allows us to apply \cite[Prop.\ 3.9]{JaMor19-1} obtaining that the identity maps are interpretations of $\vdash_{1}$ and $\vdash_{2}$ into $\vdash$. Moreover, $\vdash$ is not truth-minimal, as witnessed by the fact that $\class{K} \subseteq \ModS(\vdash)$ and $\emptyset \ne G \subsetneq F$.
\end{proof}

The proof strategy described above can be adapted to the case of truth-equational and assertional logics, as we proceed to explain.

\subsection{Truth-equational logics}

To prove that the Leibniz class of truth-equational logics is meet-prime among logics with theorems, it is convenient to recall the following characterization result.

\begin{Theorem}\label{Thm:truth-Raftery}
A logic $\vdash$ is truth-equational if and only if for every algebra $\A$,
\[
\text{if }\langle \A, F \rangle, \langle \A, G \rangle \in \ModS(\vdash) \text{, then }F = G.
\]
\end{Theorem}

\begin{proof}
See \cite[Prop.\ 17 and Thm.\ 2.8]{Ra06a}.
\end{proof}

As a consequence we obtain the desired result.

\begin{Theorem}\label{Thm:truth-eq-prime}
The Leibniz class of truth-equational logics is meet-prime among logics with theorems.
\end{Theorem}

\begin{proof}
Consider two logics $\vdash_{1}$ and $\vdash_{2}$ with theorems that are not truth-equational.  It will be enough to construct a logic $\vdash$ that is not truth-equational and in which $\vdash_{1}$ and $\vdash_{2}$ are interpretable.

First we claim that for every $i = 1, 2$ there are $\langle \A_{i}, F_{i}\rangle, \langle \A_{i}, G_{i}\rangle \in \ModS(\vdash_{i})$ such that $\emptyset \ne G_{i} \subsetneq F_{i}$.  To prove this, consider $i = 1, 2$. Since the logic $\vdash_{i}$ is not truth-equational, we can apply Theorem \ref{Thm:truth-Raftery} obtaining an algebra $\A_{i}$ and distinct $U_{i}, V_{i} \subseteq A_{i}$ such that $\langle \A_{i}, U_{i}\rangle, \langle \A_{i}, V_{i}\rangle \in \ModS(\vdash_{i})$. We can assume without loss of generality  that $V_{i} \smallsetminus U_{i} \ne \emptyset$. Bearing this in mind, we define $F_{i} \coloneqq V_{i}$ and $G_{i} \coloneqq V_{i} \cap U_{i}$. It is clear that $\langle \A_{i}, G_{i}\rangle \in \Mod(\vdash_{i})$ and $G_{i} \subsetneq F_{i}$. Since $\vdash_{i}$ has theorems, the fact that $\langle \A_{i}, G_{i}\rangle \in \Mod(\vdash_{i})$ implies that $G_{i} \ne \emptyset$. Moreover, since $G_{i} \subseteq F_{i} = V_{i}$, we have 
\[
\tarski_{\vdash_{i}}^{\A_{i}}G_{i} \subseteq \tarski_{\vdash_{i}}^{\A_{i}}V_{i}.
\]
Recall that $\tarski_{\vdash_{i}}^{\A_{i}}V_{i}$ is the identity relation, since $\langle \A_{i}, V_{i}\rangle \in \ModS(\vdash_{i})$. Hence $\tarski_{\vdash_{i}}^{\A_{i}}G_{i}$ is also the identity relation and, therefore, $\langle \A_{i}, G_{i}\rangle \in \ModS(\vdash_{i})$. Finally, by the assumptions, $\langle \A_{i}, F_{i}\rangle = \langle \A_{i}, V_{i}\rangle \in \ModS(\vdash_{i})$. This concludes the proof of the claim.

Together with Lemma \ref{Lem:truth-prime2}, the claim implies   that there are
\[
\langle \B_{1}, F \rangle, \langle \B_{1}, G \rangle \in \ModS(\vdash_{1}) \text{ and }\langle \B_{2}, F \rangle, \langle \B_{2}, G \rangle \in \ModS(\vdash_{2})
\]
for some sets $\emptyset \ne G \subsetneq F \subsetneq B$, an $\LL_{\vdash_{1}}$-algebra $\B_{1}$, and an $\LL_{\vdash_{2}}$-algebra $\B_{2}$ such that $B_{1} = B_{2} = B$.

Now, let $\B$ be the common expansion of $\B_{1}$ and $\B_{2}$ with all finitary operations on $B$. Moreover, let $\vdash$ be the logic, formulated in $\kappa \coloneqq \max\{ \vert Fm(\vdash_{1}) \vert, \vert Fm(\vdash_{2}) \vert \}$ variables induced by the set of matrices 
\[
\class{K} \coloneqq \{ \langle \B, F \rangle, \langle \B, G \rangle\}. 
\]
As in the proof of Theorem \ref{Thm:truth-min-prime}, we obtain  $\class{K} \subseteq \ModS(\vdash)$ and ${\vdash_{1}}, {\vdash_{2}} \leq {\vdash}$. Finally, applying Theorem \ref{Thm:truth-Raftery} to the fact that $\class{K} \subseteq \ModS(\vdash)$ and $F \ne G$, we conclude that $\vdash$ is not truth-equational.
\end{proof}

\subsection{Assertional logics}

\begin{Theorem}\label{Thm:asser-prime}
The Leibniz class of assertional logics is meet-prime among logics with theorems.
\end{Theorem}

\begin{proof}
Consider two logics $\vdash_{1}$ and $\vdash_{2}$ with theorems that are not assertional. As usual, it will be enough to construct a logic $\vdash$ that is not assertional and in which $\vdash_{1}$ and $\vdash_{2}$ are interpretable.

We claim that for every $i = 1, 2$ there are $\langle \A_{i}, F_{i}\rangle \in \ModS(\vdash_{i})$ such that $\vert F_{i} \vert \geq 2$ and $F_{i} \subsetneq A_{i}$. To prove this, consider $i = 1, 2$, and observe that the fact that $\vdash_{i}$ is not assertional implies that there is $\langle \A_{i}, F_{i}\rangle \in \ModS(\vdash_{i})$ such that $F_{i}$ is not a singleton. Since $\vdash_{i}$ has theorems and $\langle \A_{i}, F_{i}\rangle$ is a model of $\vdash_{i}$ we know that $F_{i} \ne \emptyset$ and, therefore, $\vert F_{i} \vert \geq 2$. This fact, together with Lemma \ref{Lem:truth-prime0}, ensures that $F_{i} \subsetneq A_{i}$, establishing the claim.

By the claim we can find a cardinal $\kappa$ large enough to guarantee that
\[
\vert A_{1}^{\kappa} \vert = \vert F_{1}^{\kappa} \vert = \vert A_{1}^{\kappa} \smallsetminus F_{1}^{\kappa} \vert = \vert A_{2}^{\kappa} \vert = \vert F_{2}^{\kappa} \vert = \vert A_{2}^{\kappa} \smallsetminus F_{2}^{\kappa} \vert.
\]
Thus there is a set $B$, a set $F \subseteq B$, an $\LL_{\vdash_{1}}$-algebra $\B_{1}$, and an $\LL_{\vdash_{2}}$-algebra $\B_{2}$ such that $B_{1} = B_{2} = B$ and
\[
\langle \B_{1}, F\rangle \cong \langle \A_{1}^{\kappa}, F_{1}^{\kappa}\rangle \text{ and }\langle \B_{2}, F\rangle \cong \langle \A_{2}^{\kappa}, F_{2}^{\kappa}\rangle.
\]
Since $\ModS(\vdash_{i})$ is closed under the formation of direct powers for every $i = 1, 2$ \cite[Lem.\ 2.3]{JaMor19-1}, we obtain $\langle \B_{1}, F\rangle \in \ModS(\vdash_{1})$ and $\langle \B_{2}, F\rangle \in \ModS(\vdash_{2})$.

Now, let $\B$ be the common expansion of $\B_{1}$ and $\B_{2}$ with all finitary operations on $B$. Moreover, let $\vdash$ be the logic, formulated in $\kappa \coloneqq \max\{ \vert Fm(\vdash_{1}) \vert, \vert Fm(\vdash_{2}) \vert \}$ variables induced by the matrix $\langle \B, F \rangle$. As in the proof of Theorem \ref{Thm:truth-min-prime}, we obtain  $\langle \B, F \rangle \in \ModS(\vdash)$ and ${\vdash_{1}}, {\vdash_{2}} \leq {\vdash}$. Finally, since $\langle \B, F \rangle$ and $F$ is not a singleton, we conclude that $\vdash$ is not assertional.
\end{proof}

\begin{Remark}
In Theorems \ref{Thm:truth-eq-prime} and \ref{Thm:asser-prime}  it is shown that the Leibniz classes of truth-equational and assertional logics are meet-prime \textit{among logics with theorems}. As we mentioned, in the light of Proposition \ref{Prop:dec-among-thm} this restriction cannot be dropped.
\qed
\end{Remark}

\subsection{Logics with theorems}

\begin{Theorem}\label{Thm:theor-prime}
The Leibniz class of logics with theorems is meet-prime.
\end{Theorem}

\begin{proof}
Let $\vdash_{\class{AI}}$ be an almost inconsistent logic, and recall that if $\vdash$ is a logic without theorems, then ${\vdash} \leq {\vdash_{\class{AI}}}$ \cite[Thm.\ 7.3]{JaMor19-1}. Together with the fact that $\vdash_{\class{IA}}$ lacks theorems, this immediately implies that the Leibniz class of logics with theorems is meet-prime.
\end{proof}

\section{Meet-reducible Leibniz classes}

Recall that protoalgebraic, equivalential, and order algebraizable logics have theorems. In the light of Proposition \ref{Prop:dec-among-thm} the corresponding Leibniz classes are trivially meet-reducible. In this section we strengthen this result by proving that they remain meet-reducible even among the restricted setting of logics with theorems. In addition, it is shown that the Leibniz class of parametrically truth-equational logics is meet-reducible in the absolute sense.\footnote{These negative results are compensated by some positive ones in Section \ref{Sec:meet-prime-logics} at least for what concerns protoalgebraic and equivalential logics.} Remarkably, the latter result cannot be inferred directly from Proposition \ref{Prop:dec-among-thm}, since parametrically truth-equational logics need not have theorems \cite[Sec.\ 4]{TMo15}. In what follows we rely on the next technical observation:

\begin{Proposition}\label{Prop:non-indexed-products}
Let $\vdash_{1}, \vdash_{2}$, and $\vdash$ be logics, and $\btau$ an interpretation of $\vdash_{1} \bigotimes \vdash_{2}$ into $\vdash$. Then for every $\langle \A, F \rangle \in \ModS(\vdash)$ there are $\langle \A_{1}, F_{1}\rangle \in \ModS(\vdash_{1})$, $\langle \A_{2}, F_{2}\rangle \in \ModS(\vdash_{2})$, and an isomorphism
\[
f \colon \langle \A_{1} \bigotimes \A_{2}, F_{1} \times F_{2}\rangle \to \langle \A^{\btau}, F \rangle.
\]
Moreover, for every submatrix $\langle \B, F \cap B\rangle \subseteq \langle \A, F \rangle$ and $i = 1, 2$, there is a submatrix $\langle \B_{i}, F_{i} \cap B_{i}\rangle \subseteq \langle \A_{i}, F_{i}\rangle$ such that:
\benroman
\item the following restriction of $f$ is a well-defined isomorphism:
\[
f \colon \langle \B_{1} \bigotimes \B_{2}, (F_{1} \cap B_{1}) \times (F_{2} \cap B_{2})\rangle \to \langle \B^{\btau}, F \cap B \rangle;
\]
\item if $F \cap B \ne \emptyset$, then for every $\theta \in \Con \B$ compatible with $F \cap B$ and every $i = 1, 2$, there is $\theta_{i} \in \Con\B_{i}$ compatible with $F_{i} \cap B_{i}$ such that 
\[
\theta = \{ \langle f\langle a, b \rangle, f \langle c, d \rangle\rangle \colon \langle a, c \rangle \in \theta_{1} \text{ and }\langle b, d \rangle \in \theta_{2} \}.
\]
\eroman
\end{Proposition}

\begin{proof}
The fact that there are $\langle \A_{1}, F_{1}\rangle \in \ModS(\vdash_{1})$, $\langle \A_{2}, F_{2}\rangle \in \ModS(\vdash_{2})$, and an isomorphism $f \colon \langle \A_{1} \bigotimes \A_{2}, F_{1} \times F_{2}\rangle \to \langle \A^{\btau}, F \rangle$ is a direct consequence of \cite[Cor.\ 4.14]{JaMor19-1}. From now on we assume without loss of generality  that $f$ is the identity map and, therefore, that
\[
\langle \A^{\btau}, F \rangle = \langle \A_{1} \bigotimes \A_{2}, F_{1} \times F_{2}\rangle.
\]

(i): Consider a submatrix $\langle \B, F \cap B\rangle \subseteq \langle \A, F \rangle$. Therefore, $\B^{\btau} \subseteq \A^{\btau} = \A_{1} \bigotimes \A_{2}$. As a consequence, for every $i = 1, 2$ there is $\B_{i} \subseteq \A_{i}$ such that $\B^{\btau} = \B_{1} \bigotimes \B_{2}$ \cite[Lem.\ 1.10]{Tay73}. Moreover, since $B = B_{1} \times B_{2}$ and $F = F_{1} \times F_{2}$, 
\[
F \cap B = (F_{1} \times F_{2}) \cap (B_{1} \times B_{2}) = (F_{1} \cap B_{1}) \times (F_{2} \cap B_{2}).
\]
As a consequence, we obtain
\[
\langle \B^{\btau}, F \cap B \rangle = \langle \B_{1} \bigotimes \B_{2}, (F_{1} \cap B_{1}) \times (F_{2} \cap B_{2})\rangle.
\]

(ii): Consider $\theta \in \Con \B$ compatible with $F \cap B$. Clearly $\theta \in \Con \B^{\btau} = \Con (\B_{1} \bigotimes \B_{2})$. As shown in \cite[Lem.\ 1.12]{Tay73}, for every $i = 1, 2$ there is $\theta_{i} \in \Con\B_{i}$ such that
\[
\theta = \{ \langle \langle a, b \rangle,  \langle c, d \rangle\rangle \colon \langle a, c \rangle \in \theta_{1} \text{ and }\langle b, d \rangle \in \theta_{2} \}.
\]
We turn to prove that $\theta_{1}$ is compatible with $F_{1} \cap B_{1}$. To this end, consider $a, c \in B_{1}$ such that $a \in F_{1} \cap B_{1}$ and $\langle a, c \rangle \in \theta_{1}$. From the assumption we have $(F_{1} \cap B_{1}) \times (F_{2} \cap B_{2}) = F \cap B \ne \emptyset$. Then there is $b \in F_{2} \cap B_{2}$. From the above display it follows that $\langle \langle a, b \rangle, \langle c, b \rangle \rangle \in \theta$. Since $\theta$ is compatible with $F$ and $\langle a, b \rangle \in F_{1} \times F_{2} = F$, we get $\langle c, b \rangle \in F = F_{1} \times F_{2}$. In particular, this guarantees that $c \in F_{1} \cap B_{1}$. As a consequence we conclude that $\theta_{1}$ is compatible with $F_{1} \cap B_{1}$. A similar argument shows that $\theta_{2}$ is compatible with $F_{2} \cap B_{2}$.
\end{proof}

\subsection{Protoalgebraic logics}

Our aim is to show that the Leibniz class of protoalgebraic logics is meet-reducible among logics with theorems. To this end, it is useful to recall a few concepts. An algebra $\A$ is said to be \textit{congruence uniform} \cite[Sec.\ 7.1]{Be11g} if $\vert a / \theta \vert = \vert b / \theta \vert$, for every $a, b \in A$ and $\theta \in \Con\A$. It is well-known that Boolean algebras are congruence uniform.

We denote by $\class{BA}$ the variety of Boolean algebras, and by $\vdash_{\class{BA}}^{\ast}$ the logic formulated in countably many variables induced by the following class of matrices:
\[
\{ \langle \A, F \rangle \colon \A \in \class{BA} \text{ and }1 \in F \}.
\]
\begin{Lemma}\label{Lem:BA*}
The logic $\vdash_{\class{BA}}^{\ast}$ has theorems, but is not protoalgebraic. Moreover, the algebraic reducts of the matrices in $\ModS(\vdash_{\class{BA}}^{\ast})$ belong to $\class{BA}$.
\end{Lemma}

\begin{proof}
Clearly $1$ is a theorem of $\vdash_{\class{BA}}^{\ast}$. Moreover, the fact that the algebraic reducts of the matrices in $\ModS(\vdash_{\class{BA}}^{\ast})$ belong to $\class{BA}$ is an immediate consequence of \cite[Cor.\ 2.6]{JaMor19-1}.

It only remains to prove that $\vdash_{\class{BA}}^{\ast}$ is not protoalgebraic. Suppose the contrary, with a view to contradiction. Then there is a set $\Delta(x, y, \vec{z})$ of congruence formulas with parameters for $\vdash_{\class{BA}}^{\ast}$. We consider the four-element Boolean algebra $\A$ with universe $\{ a, b, 0, 1 \}$, where $0$ and $1$ are respectively the bottom and the top element of the lattice order. Then we set $F \coloneqq \{ 1, a \}$ and $G \coloneqq \{ 1, a, b \}$. From the definition of $\vdash_{\class{BA}}^{\ast}$ it follows that $\langle \A, F \rangle, \langle \A, G \rangle \in \Mod(\vdash_{\class{BA}}^{\ast})$. Together with the fact that $\Delta$ is a set of congruence formulas with parameters for $\vdash_{\class{BA}}^{\ast}$ and that $F \subseteq G$, this yields that for every $p, q \in A$,
\begin{align*}
\langle p, q \rangle \in \leibniz^{\A}F \Longleftrightarrow & \, \Delta^{\A}(p, q, \vec{c}) \subseteq F \text{, for every }\vec{c}\in A\\
\Longrightarrow & \, \Delta^{\A}(p, q, \vec{c}) \subseteq G \text{, for every }\vec{c}\in A\\
\Longleftrightarrow & \, \langle p, q \rangle \in \leibniz^{\A}G.
\end{align*}
Hence we conclude that $\leibniz^{\A}F \subseteq \leibniz^{\A}G$. On the other hand, it is easy to see that $\leibniz^{\A}G$ is the identity relation, while $\leibniz^{\A}F$ is the congruence with blocks $\{ 1, a \}$ and $\{ 0, b \}$. But this contradicts the fact that $\leibniz^{\A}F \subseteq \leibniz^{\A}G$, as desired.
\end{proof}

We denote by $\class{Proto}$ the Leibniz class of protoalgebraic logics, and by $\class{Asrt}$ that of assertional logics. Bearing this in mind, we obtain the following result, in which suprema are taken  in the Leibniz hierarchy.

\begin{Theorem}\label{Thm:dec-proto}
The Leibniz class of protoalgebraic logics is meet-reducible among logics with theorems, and can be described as follows:
\[
(\class{Proto} \lor \class{Asrt}) \cap (\class{Proto} \lor \mathrm{Log}(\vdash_{\class{BA}}^{\ast})).
\]
\end{Theorem}

\begin{proof}
We set
\[
\oper{K}_{1} \coloneqq \class{Proto} \lor \class{Asrt} \;\;  \text{ and } \;\;  \oper{K}_{2} \coloneqq \class{Proto} \lor\mathrm{Log}(\vdash_{\class{BA}}^{\ast}).
\]
Observe that $\oper{K}_{1}$ and $\oper{K}_{2}$ are Leibniz classes by Theorems \ref{Thm:proto-class} and \ref{Thm:StrongAssertional}. Moreover, they comprise only logics with theorems by Propositions \ref{Prop:meets-and-joins} and \ref{Prop:still-have-thms}. We have
\begin{equation}\label{Eq:proto-1}
\class{Proto} \subsetneq \oper{K}_{1} \;\;  \text{ and }\;\;  \class{Equiv} \subsetneq \oper{K}_{2}.
\end{equation}
The validity of the inclusion $\class{Proto} \subseteq \oper{K}_{1} \cap \oper{K}_{2}$ is straightforward.\ The fact that $\class{Proto} \subsetneq \oper{K}_{1}$ is witnessed by the existence of assertional logics that are not protoalgebraic \cite[Ex.\ 7]{Ra06a}. Finally, the fact that the inclusion $\class{Equiv} \subseteq \oper{K}_{2}$ is strict follows from Lemma \ref{Lem:BA*}.

In the light of (\ref{Eq:proto-1}), it only remains to prove that $\oper{K}_{1} \cap \oper{K}_{2} \subseteq \class{Proto}$. Suppose the contrary, with a view to contradiction. Then there is a logic ${\vdash} \in \oper{K}_{1} \cap \oper{K}_{2}$ that is not protoalgebraic. By Proposition \ref{Prop:meets-and-joins}, and considering that the non-indexed product of two protoalgebraic logics is protoalgebraic and interpretable in each of them,  there are a protoalgebraic logic $\vdash_{pr}$ and an assertional logic $\vdash_{as}$ such that 
\begin{equation}\label{Eq:proto-dec2}
{\vdash_{pr} \bigotimes \vdash_{as}} \leq {\vdash} \;\; \text{ and }\:\; {\vdash_{pr} \bigotimes \vdash_{\class{BA}}^{\ast}} \leq {\vdash}.
\end{equation}
Since $\vdash$ has theorems and is not protoalgebraic, we can apply Theorem \ref{Thm:Protoalgebraic-characterization}(ii) obtaining that $\ModS(\vdash) \ne \RRR(\ModS(\vdash))$, i.e.\ that there is a matrix $\langle \A, F\rangle \in \ModS(\vdash)$ such that $\leibniz^{\A}F$ is not the identity relation.

We claim that for every pair of different $a, c \in A$, if $\langle a, c \rangle \in \leibniz^{\A}F$, then $a, c \notin F$. To prove this, consider different $a, c \in A$ such that $\langle a, c \rangle \in \leibniz^{\A}F$. By (\ref{Eq:proto-dec2}) there is an interpretation $\btau$ of $\vdash_{pr} \bigotimes \vdash_{as}$ into $\vdash$. By Proposition \ref{Prop:non-indexed-products} we obtain  without loss of generality that\footnote{For the sake of simplicity we assume that the map $f$ in the statement of Proposition \ref{Prop:non-indexed-products} is the identity. This assumption will be used systematically in this section without further notice.}
\[
\langle \A^{\btau}, F \rangle = \langle \A_{1} \bigotimes \A_{2}, F_{1} \times F_{2}\rangle
\]
for some $\langle \A_{1}, F_{1}\rangle \in \ModS(\vdash_{pr})$ and $\langle \A_{2}, F_{2}\rangle \in \ModS(\vdash_{as})$. Since $a, c \in A = A_{1} \times A_{2}$, there are $a_{1}, c_{1} \in A_{1}$ and $a_{2}, c_{2} \in A_{2}$ such that
\[
a = \langle a_{1}, a_{2}\rangle \text{ and }c = \langle c_{1}, c_{2}\rangle.
\]

From Proposition \ref{Prop:non-indexed-products}(ii) it follows that there are $\theta_{1} \in \Con \A_{1}$ and $\theta_{2} \in \Con \A_{2}$ compatible with $F_{1}$ and $F_{2}$, respectively, such that
\begin{equation}\label{Eq:proto-dec2/3}
\leibniz^{\A}F = \{ \langle \langle p, q \rangle, \langle r, s \rangle \rangle \colon \langle p, r \rangle \in \theta_{1} \text{ and }\langle q, s \rangle \in  \theta_{2}\}.
\end{equation}
Since $\vdash_{pr}$ is protoalgebraic and $\langle \A_{1}, F_{1} \rangle \in \ModS(\vdash_{pr})$, we can apply Theorem \ref{Thm:Protoalgebraic-characterization}(ii) obtaining that $\leibniz^{\A_{1}}F_{1}$ is the identity relation. Together with the fact that $\theta_{1} \subseteq \leibniz^{\A_{1}}F_{1}$ and (\ref{Eq:proto-dec2/3}), this yields 
\begin{equation}
\leibniz^{\A}F = \{ \langle \langle p, q \rangle, \langle p, s \rangle \rangle \colon p \in A_{1} \text{ and }\langle q, s \rangle \in  \theta_{2}\}.
\end{equation}
Since $\langle a, c \rangle \in \leibniz^{\A}F$ and $a \ne c$, we conclude that $a_{2} \ne c_{2}$ and $\langle a_{2}, c_{2} \rangle \in \theta_{2}$.

Now, from the fact that $\theta_{2}$ is compatible with $F_{2}$ and $\langle a_{2}, c_{2}\rangle \in \theta_{2}$, it follows  
\begin{equation}\label{Eq:proto-dec4}
\text{either }a_{2}, c_{2} \in F_{2} \text{ or }a_{2}, c_{2} \notin F_{2}.
\end{equation}
Since $\vdash_{as}$ is assertional and $\langle \A_{2}, F_{2}\rangle \in \ModS(\vdash_{as})$, we know that $F_{2}$ is a singleton. Together with $a_{2} \ne c_{2}$ and (\ref{Eq:proto-dec4}), this yields that $a_{2}, c_{2} \notin F_{2}$.  As a consequence, we obtain  $\langle a, c \rangle \notin F_{1} \times F_{2} = F$, establishing the claim.

Now, by (\ref{Eq:proto-dec2}) there is an interpretation $\btau$ of $\vdash_{pr} \bigotimes \vdash_{\class{BA}}^{\ast}$ into $\vdash$. By Proposition \ref{Prop:non-indexed-products} we obtain without loss of generality  that
\begin{align*}
\langle \A^{\btau}, F \rangle = \langle \A_{1} \bigotimes \A_{2}, F_{1} \times F_{2}\rangle
\end{align*}
for some $\langle \A_{1}, F_{1}\rangle \in \ModS(\vdash_{pr})$ and $\langle \A_{2}, F_{2}\rangle \in \ModS(\vdash_{\class{BA}}^{\ast})$.

Recall that the matrix $\langle \A, F \rangle$ is not reduced. Then there are different $a, c \in A$ such that $\langle a, c\rangle \in \leibniz^{\A}F$. Since $A = A_{1} \times A_{2}$, there are $a_{1}, c_{1} \in A_{1}$ and $a_{2}, c_{2} \in A_{2}$ such that
\[
a = \langle a_{1}, a_{2}\rangle \text{ and }c = \langle c_{1}, c_{2}\rangle.
\]
As in the proof of the claim, we obtain  
\begin{equation}\label{Eq:proto-dec5}
\leibniz^{\A}F = \{ \langle \langle p, q \rangle, \langle p, s \rangle \rangle \colon p \in A_{1} \text{ and }\langle q, s \rangle \in  \theta\}
\end{equation}
for some $\theta \in \Con\A_{2}$ compatible with $F_{2}$. Then we choose elements $e_{1} \in F_{1}$ and $e_{2} \in F_{2}$ (the fact that the logics $\vdash_{pr}$ and $\vdash_{\class{BA}}^{\ast}$ have theorems guarantees $F_{1}, F_{2} \ne \emptyset$). Clearly we have
\begin{equation}\label{Eq:proto-dec6}
\langle e_{1}, e_{2}\rangle \in F_{1} \times F_{2} = F.
\end{equation}
Together with $a \ne c$ and $\langle a, c \rangle \in \leibniz^{\A}F$, the display (\ref{Eq:proto-dec5}) implies  $a_{2} \ne c_{2}$ and $\langle a_{2}, c_{2}\rangle \in \theta$. Since $\langle \A_{2}, F_{2}\rangle \in \ModS(\vdash_{\class{BA}}^{\ast})$, we can apply Lemma \ref{Lem:BA*} obtaining that $\A_{2}$ is a Boolean algebra. Since $\A_{2}$ is congruence uniform and $\theta$ identifies two distinct elements of $A_{2}$ (namely $a_{2}$ and $c_{2}$), then there is $b \in A_{2} \smallsetminus \{ e_{2} \}$ such that $\langle e_{2}, b \rangle \in \theta$. By (\ref{Eq:proto-dec5}) we conclude  
\[
\langle \langle e_{1}, e_{2}\rangle, \langle e_{1}, b \rangle \rangle \in \leibniz^{\A}F \text{ and } \langle e_{1}, e_{2}\rangle \ne \langle e_{1}, b \rangle.
\]
Together with (\ref{Eq:proto-dec6}) this contradicts the claim. This produces the desired contradiction.
\end{proof}

\subsection{Equivalential logics}

To prove that the Leibniz class of equivalential logics is meet-reducible among logics with theorems, we need to introduce a new concept:

\begin{law}
A formula $\varphi(x)$ is an \textit{injective theorem} of a logic $\vdash$ if $\varphi(x)$ is a theorem of $\vdash$ and for every $\langle \A, F \rangle \in \ModS(\vdash)$ the term-function $\varphi^{\A} \colon A \to A$ is injective.
\end{law}

We will rely on the following:
\begin{Proposition}\label{Prop:Injective-theorem}
Logics with an injective theorem form a Leibniz class, comprising a protoalgebraic non-equivalential logic.
\end{Proposition}

The proof of the above result proceeds through a series of technical observations.

\begin{Fact}
Logics with an injective theorem form a Leibniz class.
\end{Fact}
\begin{proof}
It is not hard to show that logics with an injective theorem are closed under term-equivalence, compatible expansions, and non-indexed products of sets. In the light of Theorem  \ref{Thm:Main} this implies that they form a Leibniz class. 
\end{proof}

Consider the logic $\vdash_{\nabla}$, formulated in countably many variables and in the language consisting of a single binary connective $\to$, axiomatized by the following rules:
\[
\emptyset \rhd x \to x \qquad x, x \to y \rhd y.
\]
The logic $\vdash_{\nabla}$ has been studied in depth in \cite{JMF13ailt,Font13,JMF14opL}. We set $\nabla(x, y)\coloneqq \{ x \to y \}$.

\begin{Fact}\label{Fact:nabla-proto}
The logic $\vdash_{\nabla}$ is protoalgebraic with set of congruence formulas with parameters
\[
\hat{\nabla}(x, y, \vec{z}) \coloneqq \{ \varphi(x, \vec{z}) \to \varphi(y, \vec{z}) \colon \varphi(x, \vec{z}) \in Fm(\vdash) \}.
\]
\end{Fact}

\begin{proof}
The logic $\vdash_{\nabla}$ is protoalgebraic, since the set $\nabla(x, y)$ satisfies the requirements of Theorem \ref{Thm:Protoalgebraic-characterization}(iv). By the same theorem, $\hat{\nabla}$ is a set of congruence formulas for $\vdash_{\nabla}$.
\end{proof}

\begin{Fact}\label{Fact:implication-theorem}
For every formula $\varphi \in Fm(\vdash_{\nabla})$,
\[
\emptyset \vdash_{\nabla}\varphi \Longleftrightarrow \varphi = \psi \to \psi\text{, for some formula }\psi.
\]
\end{Fact}

\begin{proof}
See \cite[Prop.\ 2.1]{Font13}.
\end{proof}

Now, consider the extension $\vdash_{\Delta}$ of $\vdash_{\nabla}$, obtained adding for each $\psi \in \hat{\nabla}( x, y, \vec{z})$ the rule
\begin{equation}\label{Eq:admissible}
\hat{\nabla} ( x \to x, y \to y, \vec{z}) \rhd \psi. 
\end{equation}

\begin{Fact}\label{Fact:delta-proto}
The logic $\vdash_{\Delta}$ is protoalgebraic.
\end{Fact}
\begin{proof}
This is a consequence of Fact \ref{Fact:nabla-proto}, together with the fact that protoalgebraicity is preserved by extensions.
\end{proof}

\begin{Fact}
The formula $x \to x$ is an injective theorem of $\vdash_{\Delta}$.
\end{Fact}
\begin{proof}
Clearly $x \to x$ is a theorem of $\vdash_{\Delta}$. Then consider $\langle \A, F \rangle \in \ModS(\vdash_{\Delta})$ and $a, b \in A$ such that $a \to^{\A} a = b\to^{\A}b$. We need to show that $a = b$. To this end, recall by Fact \ref{Fact:delta-proto} that $\vdash_{\Delta}$ is protoalgebraic. Then we can apply Theorem \ref{Thm:Protoalgebraic-characterization}(iii), obtaining that the matrix $\langle \A, F \rangle$ is reduced, i.e.\ that the congruence $\leibniz^{\A}F$ is the identity relation. Therefore, to conclude the proof, it will be enough to show that $\langle a, b \rangle \in \leibniz^{\A}F$. By \cite[Prop.\ 2.2(i)]{JaMor19-1} this amounts to establishing that for every formula $\varphi(x, \vec{z}) \in Fm(\vdash_{\Delta})$ and every $\vec{c} \in A$,
\begin{equation}\label{Eq:admissible1}
\varphi^{\A}(a, \vec{c}) \in F \Longleftrightarrow \varphi^{\A}(b, \vec{c}) \in F.
\end{equation}

To prove the implication from left to right, consider $\varphi(x, \vec{z}) \in Fm(\vdash_{\Delta})$ and $\vec{c} \in A$ such that $\varphi^{\A}(a, \vec{c}) \in F$. Consider also an arbitrary formula $\psi(x, \vec{z}) \in Fm(\vdash_{\nabla})$. We have that $\psi(x, \vec{z}) \to \psi(x, \vec{z})$ is a theorem of $\vdash_{\Delta}$. Together with the fact that $a \to^{\A} a = b\to^{\A}b$, this yields
\[
\psi^{\A}(a \to a, \vec{c}) \to \psi^{\A}(b \to b, \vec{c}) = \psi^{\A}(a \to a, \vec{c}) \to \psi^{\A}(a \to a, \vec{c}) \in F.
\]
Hence we conclude that $\hat{\nabla}(a \to^{\A} a, b \to^{\A}b, \vec{c}) \subseteq F$. Since $\langle \A, F \rangle$ is a model of (\ref{Eq:admissible}), this implies  $\hat{\nabla}(a, b, \vec{c}) \subseteq F$. As a consequence, we obtain $\varphi^{\A}(a, \vec{c}) \to \varphi^{\A}(b, \vec{c}) \in F$. Together with the assumption that $\varphi^{\A}(a, \vec{c}) \in F$, and the fact that $\langle \A, F \rangle$ is a model of the rule $x, x \to y \rhd y$, this implies  $\varphi^{\A}(b, \vec{c}) \in F$. This concludes the proof of the left to right implication in (\ref{Eq:admissible1}). The proof of the other implication is analogous. From (\ref{Eq:admissible1}) it follows that $\langle a, b \rangle \in \leibniz^{\A}F$ and, therefore, $a = b$.
\end{proof}

Now, recall that a rule $\Gamma \rhd \varphi$ is \textit{admissible} \cite{Ry97} in a logic $\vdash$, if its addition to $\vdash$ does not produce new theorems. Equivalently, this means that $\emptyset \vdash \sigma \varphi$, for every substitution $\sigma$ such that $\emptyset \vdash \sigma[\Gamma]$.

\begin{Fact}\label{Fact:admissible}
The rules in (\ref{Eq:admissible}) are admissible in $\vdash_{\nabla}$.
\end{Fact}

\begin{proof}
To prove this, consider a substitution $\sigma$ such that $\emptyset \vdash_{\nabla} \sigma[\hat{\nabla}( x \to x, y \to y, \vec{z})]$. Since  $x \to y \in \hat{\nabla}(x, y, \vec{z})$, we have  $(x \to x) \to (y \to y) \in \hat{\nabla}(x \to x, y \to y, \vec{z})$. As a consequence, we obtain  $\emptyset \vdash_{\nabla}(\sigma x \to \sigma x) \to (\sigma y \to \sigma y)$. By Fact \ref{Fact:implication-theorem} this implies  $\sigma x \to \sigma x = \sigma y \to \sigma y$ and, therefore,  $\sigma x = \sigma y$. Since $\emptyset \vdash_{\nabla} x \to x$ and $\sigma x = \sigma y$, we obtain  $\emptyset \vdash_{\nabla} \sigma \varphi(x, \vec{z}) \to \sigma \varphi(y, \vec{z})$, for every formula $\varphi(x, \vec{z}) \in Fm(\vdash_{\nabla})$. But this amounts to the fact that $\emptyset \vdash_{\nabla} \sigma[\hat{\nabla}( x, y, \vec{z})]$. Hence we conclude that the rules in (\ref{Eq:admissible}) are admissible in $\vdash_{\nabla}$.
\end{proof}

Given a formula $\varphi$, we denote by $\Var(\varphi)$ the set of variables occurring in $\varphi$.

\begin{Fact}
The logic $\vdash_{\Delta}$ is not equivalential.
\end{Fact}

\begin{proof}
Suppose, with a view to contradiction, that $\vdash_{\Delta}$ is equivalential. Then $\vdash_{\Delta}$ has a set of congruence formulas $\Delta(x, y)$.

We claim that $\emptyset \vdash_{\Delta}\varphi$, for every formula $\varphi$ such that $\Delta(x, y) \vdash_{\Delta} \varphi$ and $\Var(\varphi) \nsubseteq \{ x, y \}$. To demonstrate this, we reason by complete induction on the length of proofs in $\vdash_{\Delta}$. Consider an ordinal $\alpha$, and suppose that $\emptyset \vdash_{\Delta} \psi$, for every formula $\psi$ such that $\Var(\varphi) \nsubseteq \{ x, y \}$, and of which the exists a proof indexed by an ordinal $< \alpha$ from $\Delta(x, y)$. Now, let $\pi \coloneqq \{ \gamma_{\beta} \colon \beta < \alpha \}$ be a proof of a formula $\varphi$ such that $\Var(\varphi) \nsubseteq \{ x, y \}$ from $\Delta(x, y)$. If $\varphi$ is a substitution instance of the axiom $x \to x$, then it is a theorem and we are done. Moreover, observe that $\varphi \notin \Delta(x, y)$, since $\Var(\varphi) \nsubseteq \{ x, y \}$. Therefore, $\varphi$ is obtained by the application of one of the inference rules of $\vdash_{\Delta}$ to a proper initial segment of $\pi$.

First consider the case where $\varphi$ is obtained by an application of the rule $x, x\to y \rhd y$. Then there is some formula $\psi$ such that  $\psi$ and $\psi \to \varphi$ appear in a proper initial segment of $\pi$. Since $\Var(\varphi) \nsubseteq \{ x, y \}$, then $\Var(\psi \to \varphi) \subseteq \{ x, y \}$. Therefore we can apply the induction hypothesis obtaining  $\emptyset \vdash_{\Delta} \psi \to \varphi$. By Fact \ref{Fact:admissible} we get  $\emptyset \vdash_{\nabla} \psi \to \varphi$. Moreover, by Fact \ref{Fact:implication-theorem} this yields  $\psi = \varphi$. In particular, this implies  $\Var(\psi) \nsubseteq \{ x, y \}$. Therefore we can apply the  induction hypothesis, obtaining  $\emptyset \vdash_{\Delta}\psi$. Since $\varphi = \psi$, we conclude that $\emptyset \vdash_{\Delta}\varphi$, as desired.

Then we consider the case where $\varphi$ is obtained by an application of one of the rules in (\ref{Eq:admissible}). Then there is a substitution $\sigma$ and a formula $\psi(x, y, \vec{z}) \in  \hat{\nabla}( x, y, \vec{z})$ such that $\varphi = \sigma \psi$, and each element of $\sigma[\hat{\nabla}( x \to x, y \to y, \vec{z})]$ appears in a proper initial segment of $\pi$. Suppose, with a view to contradiction, that $\emptyset \nvdash_{\Delta} \sigma \psi$. Since $\psi \in \hat{\nabla}(x, y, \vec{z})$, there is a formula $\gamma(x, y, \vec{z}) \in Fm(\vdash_{\nabla})$ such that $\psi = \gamma(x, y, \vec{z}) \to \gamma(y, y, \vec{z})$. Together with the fact that $\emptyset \vdash_{\Delta} x \to x$ and $\emptyset \nvdash_{\Delta} \sigma \psi$, this yields  $\sigma(\gamma(x, y, \vec{z})) \ne \sigma(\gamma(y, y, \vec{z}))$. But this easily implies  $\sigma x \ne \sigma y$ and $x \in \Var(\gamma(x, y, \vec{z}))$. As a consequence we obtain $\sigma(\gamma(x \to x, y \to y, \vec{z})) \ne \sigma(\gamma(y \to y, y \to y, \vec{z}))$. Moreover, from Facts \ref{Fact:admissible} and \ref{Fact:implication-theorem}, this yields
\[
\emptyset \nvdash_{\Delta}\sigma(\gamma(x \to x, y \to y, \vec{z})) \to \sigma(\gamma(y \to y, y \to y, \vec{z})),
\]
which amounts to 
\begin{equation}\label{Eq:contradiction-delta}
\emptyset \nvdash_{\Delta} \sigma \psi(x\to x, y, \to y, \vec{z})\text{, where }\psi = \psi(x, y, \vec{z}).
\end{equation}
Now, from the fact that $\Var(\sigma \psi) = \Var(\varphi) \nsubseteq \{ x, y \}$ it follows
\[
\Var(\sigma \psi(x \to x, y \to y, \vec{z})) \nsubseteq \{ x, y \}.
\] 
The above display and the fact that $\sigma \psi(x \to x, y \to y, \vec{z}) \in \sigma[\hat{\nabla}(x \to x, y \to y, \vec{z})]$ allow us to apply the induction hypothesis, obtaining
\[
\emptyset \vdash_{\Delta}\sigma \psi(x \to x, y \to y, \vec{z}).
\] 
But this contradicts (\ref{Eq:contradiction-delta}). Hence we conclude that $\emptyset \vdash_{\Delta} \sigma \psi$, establishing the claim.

Now we move back to the main proof. First observe that
\begin{equation}\label{Eq:delta-aaa}
\Delta(x, y) \vdash_{\Delta} (x \to z) \to (y \to z).
\end{equation}
To prove this, consider $\langle \A, F \rangle \in \Mod(\vdash_{\Delta})$ and $a, b, c \in A$ such that $\Delta^{\A}(a, b) \subseteq F$. Since $\Delta$ is a set of congruence formulas for $\vdash_{\Delta}$, this yields $\langle a, b \rangle \in \leibniz^{\A}F$ and, therefore, 
\[
\langle (a \to c) \to (a \to c), (a \to c) \to (b \to c) \rangle \in \leibniz^{\A}F.
\]
Since $x \to x$ is a theorem of $\vdash_{\Delta}$, we have that $(a \to c) \to (a \to c) \in F$. As $\leibniz^{\A}F$ is compatible with $F$, this implies  $(a \to c) \to (b \to c) \in F$, establishing (\ref{Eq:delta-aaa}). From (\ref{Eq:delta-aaa}) and the claim it follows that $\emptyset \vdash_{\Delta} (x \to z) \to (y \to z)$. By Facts \ref{Fact:admissible} and \ref{Fact:implication-theorem} this yields  $x \to z = y \to z$, which is false.
\end{proof}

We also rely on the following result \cite[Thm.\ 6.73]{AAL-AIT-f}:

\begin{Theorem}\label{Thm:SSS-equivalential}
A protoalgebraic logic $\vdash$ is equivalential if and only if $\ModS(\vdash)$ is closed under \nolinebreak$\SSS$.
\end{Theorem}

We denote by $\class{Equiv}$ and $\class{Injctv}$ the Leibniz classes of equivalential logics and of logics with an an injective theorem, respectively. Bearing this in mind, the main result of this part takes the following form:

\begin{Theorem}\label{Thm:dec-equivalential}
The Leibniz class of equivalential logics is meet-reducible among logics with theorems, and can be described as follows:
\[
(\class{Equiv} \lor (\class{Proto} \cap \class{Injctv})) \cap (\class{Equiv} \lor (\class{Proto} \cap \class{Asrt})).
\]
\end{Theorem}

\begin{proof}
We begin by setting
\[
\oper{K}_{1} \coloneqq \class{Equiv} \lor (\class{Proto} \cap \class{Injctv}) \text{ and }\oper{K}_{2} \coloneqq \class{Equiv} \lor (\class{Proto} \cap \class{Asrt}).
\]
The fact that $\oper{K}_{1}$ and $\oper{K}_{2}$ are Leibniz classes is a consequence of Theorems \ref{Thm:proto-class}, \ref{Thm:equiv-class}, and \ref{Thm:StrongAssertional}, and Proposition \ref{Prop:Injective-theorem}. Moreover, they comprise only logics with theorems by Propositions \ref{Prop:meets-and-joins} and \ref{Prop:still-have-thms}. We have
\begin{equation}\label{Eq:equivalential-1}
\class{Equiv} \subsetneq \oper{K}_{1} \text{ and }\class{Equiv} \subsetneq \oper{K}_{2}.
\end{equation}
The validity of the inclusion $\class{Equiv} \subseteq \oper{K}_{1} \cap \oper{K}_{2}$ is straightforward.\ Moreover, Proposition \ref{Prop:Injective-theorem} guarantees that   $\class{Equiv} \subsetneq \oper{K}_{1}$. Finally, the fact that the inclusion $\class{Equiv} \subseteq \oper{K}_{2}$ is strict follows from the following observation \cite[Prop.\ 6.1 and Thm.\ 6.4]{CzJa00}:
\[
(\class{Proto} \cap \class{Asrt}) \smallsetminus \class{Equiv} \ne \emptyset.
\]

In the light of (\ref{Eq:equivalential-1}), it only remains to prove  $\oper{K}_{1} \cap \oper{K}_{2} \subseteq \class{Equiv}$. Suppose the contrary, with a view to contradiction. Then there is a logic ${\vdash} \in \oper{K}_{1} \cap \oper{K}_{2}$ that is not equivalential. By Proposition \ref{Prop:meets-and-joins} there are an equivalential logic $\vdash_{eq}$, a protoalgebraic logic with an injective theorem $\vdash_{pin}$, and an assertional protoalgebraic logic $\vdash_{ap}$ such that

\begin{equation}\label{Eq:equivalential-dec0}
{\vdash_{eq} \bigotimes \vdash_{pin}} \leq {\vdash} \text{ and }{\vdash_{eq} \bigotimes \vdash_{ap}} \leq {\vdash}.
\end{equation}
Together with \cite[Prop.\ 3.8]{JaMor19-1} this implies that $\vdash$ is term-equivalent to a compatible expansion of the non-indexed product of a pair of protoalgebraic logics. Since protoalgebraic logics form a Leibniz class, by Theorem \ref{Thm:Main}(ii) we conclude that $\vdash$ is protoalgebraic. Moreover, since $\vdash$ is not equivalential, we can apply Theorem \ref{Thm:SSS-equivalential} obtaining that $\ModS(\vdash)$ is not closed under $\SSS$, i.e.\ there is a matrix $\langle \A, F\rangle \in \ModS(\vdash)$ with a submatrix $\langle \B, F \cap B \rangle \subseteq \langle \A, F\rangle$ such that $\langle \B, F \cap B \rangle \notin \ModS(\vdash)$.

We claim that for every pair of different $a, c \in B$, if $\langle a, c \rangle \in \leibniz^{\B}(F \cap B)$, then $a, b \notin F$. To prove this, consider different $a, c \in B$ such that $\langle a, c \rangle \in \leibniz^{\B}(F \cap B)$. By (\ref{Eq:equivalential-dec2}) there is an interpretation $\btau$ of $\vdash_{eq} \bigotimes \vdash_{ap}$ into $\vdash$. By Proposition \ref{Prop:non-indexed-products}(i) we obtain without loss of generality that
\begin{align*}
\langle \A^{\btau}, F \rangle &= \langle \A_{1} \bigotimes \A_{2}, F_{1} \times F_{2}\rangle\\
\langle \B^{\btau}, F \cap B \rangle &= \langle \B_{1} \bigotimes \B_{2}, (F_{1} \cap B_{1}) \times (F_{2} \cap B_{2}) \rangle
\end{align*}
for some $\langle \B_{1}, F_{1}\rangle \subseteq \langle \A_{1}, F_{1}\rangle \in \ModS(\vdash_{eq})$ and $\langle \B_{2}, F_{2}\rangle \subseteq \langle \A_{2}, F_{2}\rangle \in \ModS(\vdash_{pa})$. Since $a, c \in B = B_{1} \times B_{2}$, there are $a_{1}, c_{1} \in B_{1}$ and $a_{2}, c_{2} \in B_{2}$ such that
\[
a = \langle a_{1}, a_{2}\rangle \text{ and }c = \langle c_{1}, c_{2}\rangle.
\]
Again from Proposition \ref{Prop:non-indexed-products}(ii) it follows that there are $\theta_{1} \in \Con \B_{1}$ and $\theta_{2} \in \Con \B_{2}$ compatible with $F_{1} \cap B_{1}$ and $F_{2} \cap B_{2}$, respectively, such that
\[
\leibniz^{\B}(F \cap B) = \{ \langle \langle p, q \rangle, \langle r, s \rangle \rangle \colon \langle p, r \rangle \in \theta_{1} \text{ and }\langle q, s \rangle \in  \theta_{2}\}.
\]
Together with the fact that $\langle a, c \rangle \in \leibniz^{\B}(F \cap B)$, this yields
\begin{equation}\label{Eq:equivalential-dec2}
\langle a_{i}, c_{i}\rangle \in \theta_{i} \subseteq \leibniz^{\B_{i}}(F_{i} \cap B_{i}), \text{ for every }i = 1, 2.
\end{equation}

Since $\vdash_{eq}$ is equivalential, we can apply Theorem \ref{Thm:SSS-equivalential} obtaining that $\ModS(\vdash_{eq})$ is closed under $\SSS$. Together with  $\langle \A_{1}, F_{1}\rangle \in \ModS(\vdash_{eq})$, this implies  $\langle \B_{1}, F_{1} \cap B_{1}\rangle \in \ModS(\vdash_{eq})$. Now, recall from Theorem \ref{Thm:Protoalgebraic-characterization}(ii) that $\ModS(\vdash_{eq}) = \RRR(\Mod(\vdash_{eq}))$. As a consequence, the congruence $\leibniz^{\B_{1}}(F_{1} \cap B_{1})$ is the identity relation. Hence by (\ref{Eq:equivalential-dec2}) we get  $a_{1} = c_{1}$.  Since $a \ne c$, we conclude that $a_{2} \ne c_{2}$. 

Now, from (\ref{Eq:equivalential-dec2}) and the fact that $\leibniz^{\B_{2}}(F_{2} \cap B_{2})$ is compatible with $F_{2} \cap B_{2}$ it follows that 
\begin{equation}\label{Eq:equivalential-dec3}
\text{either }a_{2}, c_{2} \in F_{2} \cap B_{2} \text{ or }a_{2}, c_{2} \notin F_{2} \cap B_{2}.
\end{equation}
Since $\vdash_{pa}$ is assertional and $\langle \A_{2}, F_{2}\rangle \in \ModS(\vdash_{pa})$, we know that $F_{2}$ is a singleton. Together with $a_{2} \ne c_{2}$ and (\ref{Eq:equivalential-dec3}), this yields  $a_{2}, c_{2} \notin F_{2} \cap B_{2}$.  As a consequence, we obtain  $\langle a, c \rangle \notin F_{1} \times F_{2} = F$, establishing the claim.

Now, by (\ref{Eq:equivalential-dec2}) there is an interpretation $\btau$ of $\vdash_{eq} \bigotimes \vdash_{pin}$ into $\vdash$. By Proposition \ref{Prop:non-indexed-products}(i) we obtain without loss of generality that
\begin{align*}
\langle \A^{\btau}, F \rangle &= \langle \A_{1} \bigotimes \A_{2}, F_{1} \times F_{2}\rangle\\
\langle \B^{\btau}, F \cap B \rangle &= \langle \B_{1} \bigotimes \B_{2}, (F_{1} \cap B_{1}) \times (F_{2} \cap B_{2}) \rangle
\end{align*}
for some $\langle \B_{1}, F_{1}\rangle \subseteq \langle \A_{1}, F_{1}\rangle \in \ModS(\vdash_{eq})$ and $\langle \B_{2}, F_{2}\rangle \subseteq \langle \A_{2}, F_{2}\rangle \in \ModS(\vdash_{pin})$.

Recall that the matrix $\langle \B, F \cap B \rangle$ is not reduced. Then there are different $a, c \in B$ such that $\langle a, c\rangle \in \leibniz^{\B}(F \cap B)$. Since $B = B_{1} \times B_{2}$, there are $a_{1}, c_{1} \in B_{1}$ and $a_{2}, c_{2} \in B_{2}$ such that
\[
a = \langle a_{1}, a_{2}\rangle \text{ and }c = \langle c_{1}, c_{2}\rangle.
\]
As in the proof of the claim, we obtain  $a_{2} \ne c_{2}$ and $\langle a_{2}, c_{2}\rangle \in \leibniz^{\B_{2}}(F_{2} \cap B_{2})$. 

Now, let $\varphi(x)$ be an injective theorem of $\vdash_{pin}$. Since $\langle \A_{2}, F_{2}\rangle \in \ModS(\vdash_{pin})$ and $a_{2} \ne c_{2}$, we have 
\begin{equation}\label{Eq:equivalential-dec4}
\varphi^{\A_{2}}(a_{2}) \ne \varphi^{\A_{1}}(a_{2}) \text{ and }\varphi^{\A_{2}}(a_{2}), \varphi^{\A_{2}}(c_{2}) \in F_{2}. 
\end{equation}
Let also $\top(x)$ be an arbitrary theorem of $\vdash_{eq}$. Since $\langle \A_{1}, F_{1}\rangle \in \Mod(\vdash_{eq})$, we have  
\begin{equation}\label{Eq:equivalential-dec5}
\top^{\A_{1}}(a_{1}), \top^{\A_{1}}(c_{1}) \in F_{1}.
\end{equation}
Observe that the pair $\langle \top, \varphi\rangle$ is a unary connective of ${\vdash_{eq}} \bigotimes {\vdash_{pld}}$. Together with (\ref{Eq:equivalential-dec4}, \ref{Eq:equivalential-dec5}), this yields 
\begin{align*}
\langle \top, \varphi \rangle^{\B_{1} \bigotimes \B_{2}} (a)&= \langle \top^{\A_{1}}(a_{1}), \varphi^{\A_{2}}(a_{2})\rangle \in F_{1} \times F_{2} = F\\
\langle \top, \varphi \rangle^{\B_{1} \bigotimes \B_{2}} (c)&= \langle \top^{\A_{1}}(c_{1}), \varphi^{\A_{2}}(c_{2})\rangle \in F_{1} \times F_{2} = F
\end{align*}
and
\[
\langle \top, \varphi \rangle^{\B_{1} \bigotimes \B_{2}} (a) \ne \langle \top, \varphi \rangle^{\B_{1} \bigotimes \B_{2}} (c).
\]
Since $\B_{1} \bigotimes \B_{2} = \B^{\btau}$, we know that $\langle \top, \varphi \rangle^{\B_{1} \bigotimes \B_{2}}$ is a term-function of $\B$. Together with the above displays and the fact that $\langle a, c \rangle \in \leibniz^{\B}(F \cap B)$, this implies that $\leibniz^{\B}(F \cap B)$ identifies two different elements of $F$, i.e.\ $\langle \top, \varphi \rangle^{\B_{1} \bigotimes \B_{2}} (a)$ and $\langle \top, \varphi \rangle^{\B_{1} \bigotimes \B_{2}} (c)$. But this contradicts the claim. Hence we reached a contradiction, as desired.
\end{proof}

\subsection{Order algebraizable logics}

We denote by $\class{Order}$ and $\class{Truth}$ the Leibniz classes of order algebraizable and truth-equational logics, respectively.

\begin{Theorem}\label{Thm:OrderAlgDecomposable}
The Leibniz class of order algebraizable logics is meet-reducible among logics with theorems, and can be described as follows:
\[
(\class{Order} \lor \class{Truth}) \cap \class{Equiv}.
\]
\end{Theorem}

\begin{proof}
First we set $\oper{K} \coloneqq \class{Order} \lor \class{Truth}$. Observe that $\oper{K}$ and $\class{Equiv}$ are Leibniz classes by Theorems \ref{Thm:equiv-class}, \ref{Thm:truth-class}, and \ref{Thm:order-class}. The fact that they comprise only logics with theorems is a consequence of Propositions \ref{Prop:meets-and-joins} and \ref{Prop:still-have-thms}. Moreover, we have
\begin{equation}\label{Eq:order-1}
\class{Order} \subsetneq \oper{K}\;\; \text{ and }\;\; \class{Order} \subsetneq \class{Equiv}.
\end{equation}
To prove this, recall that every order algebraizable logic is equivalential \cite[Prop.\ 7.11(iii)]{JRa13}. In particular, this implies that $\class{Order} \subseteq \oper{K} \cap \class{Equiv}$.\ The fact that the inclusion $\class{Order} \subseteq \oper{K}$ is strict is an immediate consequence of the fact that so is the inclusion $\class{Equiv} \subseteq \oper{K}_{2}$ in the proof of Theorem \ref{Thm:dec-equivalential}. On the other hand, the fact that $\class{Order} \subsetneq \class{Equiv}$ is witnessed by the existence of equivalential logics that are not order algebraizable \cite[p.\ 267]{JRa13}.

In the light of (\ref{Eq:order-1}), it only remains to prove that $\oper{K} \cap \class{Equiv} \subseteq \class{Order}$. To this end, consider a logic ${\vdash} \in \oper{K} \cap \class{Equiv}$. Clearly, $\vdash$ is equivalential and, therefore, there is a set of formulas $\Delta(x, y)$ of $\vdash$ such that for every $\langle \A, F \rangle \in \ModS(\vdash)$ and $a, c \in A$,
\begin{equation}\label{Eq:order-2}
a = c \Longleftrightarrow \Delta^{\A}(a, c) \subseteq F.
\end{equation}
Moreover, since ${\vdash} \in \oper{K}$, there are an order algebraizable logic $\vdash_{or}$, a truth-equational logic $\vdash_{tr}$, and an interpretation $\btau$ of $\vdash_{or} \bigotimes \vdash_{tr}$ into $\vdash$. Then there there are a set of formulas $\nabla(x, y)$ and a set of inequalities $I(x)$ of $\vdash_{or}$ such that for every $\langle \A, F \rangle \in \ModS(\vdash_{or})$,
\benroman
\item the relation $\preccurlyeq^{\A}_{F}$ on $A$ defined  for every $a, c \in A$ by
\[
a \preccurlyeq_{F}^{\A} c \Longleftrightarrow \nabla(a, c)^{\A} \subseteq F
\]
is a partial order; and
\item for every $a \in A$,  $a \in F$ if and only if $\langle \A, \preccurlyeq_{F}^{\A} \rangle \vDash I(a)$.
\eroman
Finally, since $\vdash_{tr}$ is truth-equational, there is a set of equations $E(x)$ of $\vdash_{tr}$ such that for every $\langle \A, F \rangle \in \ModS(\vdash_{tr})$ and $a \in A$
\begin{equation}\label{Eq:order-3}
a \in F \Longleftrightarrow \A \vDash E(a).
\end{equation}

Now, we choose a theorem $\top$ of $\vdash_{tr}$, and for every formula $\varphi \in \nabla(x, y)$ we consider the following basic operations of ${\vdash_{or}} \bigotimes {\vdash_{tr}}$, in which $\pi_{i}$ is the projection map on the $i$-th coordinate:
\begin{align*}
\hat{\varphi}(x, y) &\coloneqq \langle \varphi(x_{1}, x_{2}), \top(x_{1})\rangle\\
x \multimap y &\coloneqq \langle \pi_{1}(x_{1}, x_{2}), \pi_{2}(x_{1}, x_{2})\rangle.
\end{align*}
Observe that for every $\LL_{\vdash_{or}}$-algebra $\A_{1}$, $\LL_{\vdash_{tr}}$-algebra $\A_{2}$, and elements $\langle a_{1}, a_{2}\rangle, \langle c_{1}, c_{2}\rangle \in A_{1} \times A_{2}$, we have 
\begin{align}\label{Eq:order-4}
    \begin{split}
    \hat{\varphi}^{\A_{1}\bigotimes \A_{2}}(\langle a_{1}, a_{2}\rangle, \langle c_{1}, c_{2}\rangle) &= \langle \varphi^{\A_{1}}(a_{1}, c_{1}), \top^{\A_{2}}(a_{2})\rangle\\
\langle a_{1}, a_{2}\rangle \multimap^{\A_{1}\bigotimes \A_{2}} \langle c_{1}, c_{2}\rangle &= \langle a_{1}, c_{2}\rangle.
\end{split}
\end{align}
Then we define the following sets of formulas and inequalities of $\vdash$:
\begin{align*}
\nabla^{\ast}(x, y) &\coloneqq \{ \btau(\hat{\varphi}) \colon \varphi \in \nabla \} \cup \Delta(x, \btau(x \multimap y)) \\
I^{\ast}(x) & \coloneqq \{ \btau(\langle \epsilon(x_{1}), \zeta(x_{1})\rangle) \preccurlyeq \btau( \langle \delta(x_{1}), \gamma(x_{1})\rangle) \colon \epsilon \preccurlyeq \delta \in I \text{ and }\zeta \thickapprox \gamma \in E \}.
\end{align*}

To conclude the proof, it will be enough to show that the sets $\nabla^{\ast}$ and $I^{\ast}$ witness the order algebraizability of $\vdash$. To this end, consider a matrix $\langle \A, F \rangle \in \ModS(\vdash)$, and let $\preccurlyeq_{F}^{\A}$ be the relation on $A$ defined for every $a, c \in A$ as
\[
a \preccurlyeq_{F}^{\A} c \Longleftrightarrow \nabla^{\ast}(a, c)^{\A} \subseteq F.
\]
We need to show that $\preccurlyeq_{F}^{\A}$ is a partial order on $A$ and that for every $a \in A$,
\begin{equation}\label{Eq:order-5}
a \in F \Longleftrightarrow \langle \A, \preccurlyeq_{F}^{\A}\rangle \vDash I^{\ast}(a).
\end{equation}

We claim that for every $\langle a_{1}, a_{2}\rangle, \langle c_{1}, c_{2}\rangle \in A_{1} \times A_{2} = A$,
\[
\langle a_{1}, a_{2}\rangle \preccurlyeq_{F}^{\A}\langle c_{1}, c_{2}\rangle \Longleftrightarrow a_{1} \preccurlyeq_{F_{1}}^{\A_{1}} c_{1} \text{ and }a_{2} = c_{2}.
\]
To prove this, observe that
\begin{align*}
\langle a_{1}, a_{2}\rangle \preccurlyeq_{F}^{\A}\langle c_{1}, c_{2}\rangle  \Longleftrightarrow & \, \nabla^{\ast}(\langle a_{1}, a_{2}\rangle, \langle c_{1}, c_{2}\rangle)^{\A} \subseteq F\\
 \Longleftrightarrow & \, \{ \btau(\hat{\varphi})^{\A}(\langle a_{1}, a_{2}\rangle, \langle c_{1}, c_{2}\rangle) \colon \varphi \in \nabla \} \subseteq F \text{ and}\\
& \, \Delta^{\A}(\langle a_{1}, a_{2}\rangle, \btau(\langle a_{1}, a_{2}\rangle \multimap \langle c_{1}, c_{2}\rangle) \subseteq F \\
\Longleftrightarrow & \, \{ \hat{\varphi}^{\A_{1} \bigotimes \A_{2}}(\langle a_{1}, a_{2}\rangle, \langle c_{1}, c_{2}\rangle) \colon \varphi \in \nabla \} \subseteq F\text{ and}\\
& \, \Delta^{\A}(\langle a_{1}, a_{2}\rangle, \langle a_{1}, a_{2}\rangle \multimap^{\A_{1} \bigotimes \A_{2}} \langle c_{1}, c_{2}\rangle) \subseteq F\\
\Longleftrightarrow & \, \{ \langle \varphi^{\A_{1}}(a_{1}, c_{1}), \top^{\A_{2}}(a_{2})\rangle \colon \varphi \in \nabla \} \subseteq F\text{ and}\\
& \, \Delta^{\A}(\langle a_{1}, a_{2}\rangle, \langle a_{1}, c_{2}\rangle) \subseteq F\\
\Longleftrightarrow & \, \nabla^{\A_{1}}(a_{1}, c_{1}) \subseteq F_{1} \text{ and }\langle a_{1}, a_{2}\rangle = \langle a_{1}, c_{2}\rangle\\
\Longleftrightarrow & \, a_{1} \preccurlyeq_{F_{1}}^{\A_{1}} c_{1} \text{ and }a_{2} = c_{2}.
\end{align*}
The above equivalences are justified as follows: the first, the second, and the sixth are straightforward, the third follows from the fact that $\A^{\btau} = \A_{1} \bigotimes \A_{2}$, the fourth is a consequence of (\ref{Eq:order-4}), and the fifth follows from (\ref{Eq:order-2}) and the observation that $F = F_{1} \times F_{2}$ and $\top^{\A_{2}}(a_{2}) \in F_{2}$. This establishes the claim.

Recall that $\preccurlyeq_{F_{1}}^{\A_{1}}$ is a partial order on $A_{1}$ by (i). Together with the claim, this implies that $\preccurlyeq_{F}^{\A}$ is a partial order on $A$. Then consider an element $\langle a_{1}, a_{2}\rangle \in A_{1} \times A_{2} = A$. We have 
\begin{align*}
\langle a_{1}, a_{2}\rangle \in F \Longleftrightarrow & \,  a_{1} \in F_{1} \text{ and }a_{2} \in F_{2}\\
\Longleftrightarrow& \, \langle \A_{1}^{\btau}, \preccurlyeq_{F_{1}}^{\A_{1}} \rangle \vDash I(a_{1}) \text{ and }\A_{2}^{\btau} \vDash E(a_{2})\\
 \Longleftrightarrow & \, \langle \A, \preccurlyeq_{F}^{\A} \rangle \vDash I^{\ast}(\langle a_{1}, a_{2}\rangle).
\end{align*}
The above equivalences are justified as follows: the first is a consequence of the equality $F = F_{1} \times F_{2}$, the second follows from (ii) and (\ref{Eq:order-3}), and the third from the claim. This establishes (\ref{Eq:order-5}). Hence we conclude that $\vdash$ is order algebraizable.
\end{proof}

\subsection{Parametrically truth-equational logics}

As we mentioned, parametrically truth-equational logics need not have theorems in general \cite[Sec.\ 4]{TMo15}. In particular, they lie outside the scope of Proposition \ref{Prop:dec-among-thm} and we cannot immediately infer that their Leibniz class is meet-reducible in the absolute sense. We proceed to prove that this is indeed the case. To this end, we need the following observation:

\begin{Proposition} \label{Prop:mnml-inclusion}
Every parametrically truth-equational logic is truth-minimal, but the converse does not hold in general.
\end{Proposition}

\begin{proof}
Let $\vdash$ be a parametrically truth-equational logic. Then consider two matrices $\langle \A, F \rangle, \langle \A, G \rangle \in \ModS(\vdash)$ such that $\emptyset \ne G \subseteq F$. Let also $E(x, \vec{y})$ be the set of equations that witnesses the fact that $\vdash$ is parametrically truth-equational. Since $\langle \A, F \rangle, \langle \A, G \rangle \in \ModS(\vdash)$ and $F, G \ne \emptyset$, for every $a \in A$ we have 
\[
a \in F \Longleftrightarrow \A \vDash E(a, \vec{c}) \text{ for every }\vec{c} \in A \Longleftrightarrow a \in G.
\]
As a consequence we obtain  $F = G$ and, therefore, that $\vdash$ is truth-minimal. This establishes that every parametrically truth-equational logic is truth-minimal.

To conclude the proof, we need to exhibit a truth-minimal logic that is not parametrically truth-equational. To this end, let $\vdash$ be the logic formulated in countably many variables induced by the set of matrices $\{ \langle \B_{2}, \{ 1 \} \rangle, \langle \B_{2}, \{ 0 \} \rangle \}$, where $\B_{2}$ is the two-element Boolean algebra with universe $\{ 0, 1 \}$. By Proposition \ref{Prop:min-trick3} the logic $\vdash$ is truth-minimal. Now, since $\B_{2}$ is a two-element algebra, it is immediate that the matrices $\langle \B_{2}, \{ 1 \} \rangle$ and $\langle \B_{2}, \{ 0 \} \rangle$ are reduced. As a consequence, we obtain 
\[
\langle \B_{2}, \{ 1 \} \rangle, \langle \B_{2}, \{ 0 \} \rangle \in \RRR(\Mod(\vdash)) \subseteq \ModS(\vdash).
\]
Suppose, with a view to contradiction, that $\vdash$ is parametrically truth-equational, and let $E(x, \vec{y})$ be the set of equations witnessing this fact. We have 
\begin{align*}
1 \in \{ 1 \} &\Longrightarrow \B_{2} \vDash E(1, \vec{c}), \text{ for every }\vec{c} \in B_{2}\\
& \Longrightarrow 1 \in \{ 0 \}\\
& \Longrightarrow 0 = 1.
\end{align*}
The first implication above follows from the fact that $\langle \B_{2}, \{ 1 \} \rangle \in \ModS(\vdash)$, the second from $\langle \B_{2}, \{ 0 \} \rangle \in \ModS(\vdash)$, and the third is straightforward. Since $0 \ne 1$, this produces a contradiction. Hence we conclude that $\vdash$ is not parametrically truth-equational.
\end{proof}

Given a logic $\vdash$ and an $\LL_{\vdash}$-algebra $\A$, we denote by $\FFi_{\vdash}\A$ the set of deductive filters of $\vdash$ on $\A$. We build on the following  characterization result.

\begin{Theorem}\label{Thm:charac-parametrically}
A logic $\vdash$ is parametrically truth-equational if and only if for every $\LL_{\vdash}$-algebra $\A$ and every family $X \cup \{ F \} \subseteq \FFi_{\vdash}\A \smallsetminus \{ \emptyset \}$,
\[
\text{if }\bigcap\{  \leibniz^{\A}G \colon G \in X \} \subseteq \leibniz^{\A}F\text{, then }\bigcap X \subseteq F.
\] 
\end{Theorem}

\begin{proof}
The result is essentially a re-working of an analogous characterization of truth-equational logics in \cite{Ra06a}. For the details, see \cite[Thm.\ 3.9]{TMo15}. 
\end{proof}

We denote by $\class{Thrms}$, $\class{ParTruth}$, and $\class{Mnml}$ the Leibniz classes of logics with theorems, parametrically truth-equational logics, and truth-minimal logics respectively.

\begin{Theorem}\label{Thm:red-parametrically}
The Leibniz class of parametrically truth-equational logics is meet-reducible, and can be described as follows:
\[
(\class{Thrms} \lor \class{ParTruth}) \cap \class{Mnml}.
\]
\end{Theorem}

\begin{proof}
First we set $\oper{K} \coloneqq \class{Thrms} \lor \class{ParTruth}$, and observe that $\oper{K}$ and $\class{Mnml}$ are Leibniz classes by Theorems \ref{Thm:truth-class}, \ref{Thm:order-class}, and \ref{Thm:truth-minimal-class}. Moreover, we have
\begin{equation}\label{Eq:par-1}
\class{ParTruth} \subsetneq \oper{K}\;\;  \text{ and }\;\; \class{ParTruth} \subsetneq \class{Mnml}.
\end{equation}
The validity of the inclusion $\class{ParTruth} \subseteq \oper{K}$ is straightforward. The fact that it is strict is witnessed by the existence of logics with theorems that are not parametrically truth-equational, e.g.\ \cite[Ex.\ 7.5]{TMo15}. On the other hand, from Proposition \ref{Prop:mnml-inclusion} it follows that $\class{ParTruth} \subsetneq \class{Mnml}$.

In the light of (\ref{Eq:par-1}), it only remains to prove that $\oper{K} \cap \class{Mnml} \subseteq \class{ParTruth}$. To this end, consider a logic ${\vdash} \in \oper{K} \cap \class{Mnml}$. Clearly, $\vdash$ is truth-minimal, and there are a parametrically truth-equational logic $\vdash_{pt}$, and a logic with theorems $\vdash_{thm}$ such that ${\vdash_{pt} \bigotimes \vdash_{thm}} \leq {\vdash}$. Then let $\btau$ be an interpretation of $\vdash_{pt} \bigotimes \vdash_{thm}$ into $\vdash$.

In order to establish that $\vdash$ is parametrically truth-equational, we rely on Theorem \ref{Thm:charac-parametrically}. Consider an $\LL_{\vdash}$-algebra on $\A$, and a family $X \cup \{ F \} \subseteq \FFi_{\vdash} \A \smallsetminus \{ \emptyset \}$ such that 
\begin{equation}\label{Eq:par-2}
\bigcap \{ \leibniz^{\A}G \colon G \in X \} \subseteq \leibniz^{\A}F.
\end{equation}
To conclude the proof it suffices to show that $\bigcap X \subseteq F$. If $\bigcap X = \emptyset$, we are done. Then we consider the case where $\bigcap X$ is non-empty. For the sake of readability, we set $H \coloneqq \bigcap X$ and $\theta \coloneqq \tarski_{\vdash}^{\A} \bigcap X$.

From the fact that $\langle \A / \theta, H / \theta \rangle \in \ModS(\vdash)$ and \cite[Prop.\ 4.12]{AFRM15}, it follows that
\begin{equation}\label{Eq:par-3}
\langle (\A / \theta)^{\btau}, H / \theta \rangle = \langle \A_{1}, H_{1}\rangle \bigotimes \langle \A_{2}, H_{2}\rangle
\end{equation}
for some $\langle \A_{1}, H_{1}\rangle \in \ModS(\vdash_{pt})$ and $\langle \A_{2}, H_{2}\rangle \in \ModS(\vdash_{thm})$ with $H_{1}, H_{2} \ne \emptyset$. Moreover, by (\ref{Eq:par-2}) we obtain 
\[
\theta \subseteq \bigcap \{ \leibniz^{\A}G \colon G \in X \} \subseteq \leibniz^{\A}F.
\]
As a consequence, $\theta$ is compatible with $F$ and, therefore, $F /\theta \in \FFi_{\vdash}(\A / \theta)$. Together with the fact that $\btau$ is an interpretation of $\vdash_{pt} \bigotimes \vdash_{thm}$ into $\vdash$ and \cite[Prop.\ 3.3]{JaMor19-1}, this yields
\[
F  /\theta \in \FFi_{\vdash_{pt} \bigotimes \vdash_{thm}}((\A / \theta)^{\btau}) \text{ and }F / \theta \ne  \emptyset.
\]
With an application of Lemma \ref{Lem:min-trick1} to the above display and (\ref{Eq:par-3}), we conclude that
\begin{equation}\label{Eq:par-4}
F / \theta = F_{1} \times F_{2} \text{ for some }F_{1} \in \FFi_{\vdash_{pt}}\A_{1} \smallsetminus \{ \emptyset \} \text{, and }F_{2} \in \FFi_{\vdash_{thm}}\A_{2}.
\end{equation}

We claim that $H / \theta \cap F/ \theta  \ne \emptyset$. To prove this, consider the family $Y \coloneqq \{ G \in \FFi_{\vdash_{pt}}\A_{1} \colon H_{1} \subseteq G\}$. We have
\[
\bigcap \{ \leibniz^{\A_{1}}G \colon G \in Y \} = \tarski_{\vdash_{pt}}^{\A_{1}}H_{1} \subseteq  \leibniz^{\A_{1}}F_{1},
\]
where the inclusion $\tarski_{\vdash_{pt}}^{\A_{1}}H_{1} \subseteq  \leibniz^{\A_{1}}F_{1}$ follows from the fact that $\tarski_{\vdash_{pt}}^{\A_{1}}H_{1}$ is the identity relation, because $\langle \A_{1}, H_{1}\rangle \in \ModS(\vdash_{pt})$. Since $Y \cup \{ F_{1}\} \subseteq \FFi_{\vdash_{pt}}\A_{1} \smallsetminus \{ \emptyset \}$, we apply the fact that $\vdash_{pt}$ is parametrically truth-equational and Theorem \ref{Thm:charac-parametrically}, obtaining $H_{1} = \bigcap Y \subseteq F_{1}$. In particular, this guarantees that $H_{1} \cap F_{1} \ne \emptyset$, since $H_{1} \ne \emptyset$. Similarly, the fact that $\vdash_{thm}$ has theorems and $H_{2}, F_{2} \in \FFi_{\vdash_{thm}}\A_{2}$ implies $H_{2} \cap F_{2} \ne \emptyset$. Thus for every $i = 1,2$ there is $a_{i} \in H_{i} \cap F_{i}$. By (\ref{Eq:par-3}, \ref{Eq:par-4}) we conclude
\[
\langle a_{1}, a_{2}\rangle \in (H_{1} \times H_{2}) \cap (F_{1} \times F_{2}) = H/ \theta \cap F / \theta,
\]
establishing the claim.

Now, observe that the intersection $Z \coloneqq H / \theta \cap F/ \theta$ is clearly a deductive filter of $\vdash$ on $\A / \theta$. Together with the fact that $Z \subseteq H / \theta$, this implies $\tarski_{\vdash}^{\A / \theta}Z \subseteq \tarski_{\vdash}^{\A/ \theta}(H / \theta)$. Since the latter congruence is the identity, so is $\tarski_{\vdash}^{\A}Z$ and, therefore, $\langle \A / \theta, Z \rangle \in \ModS(\vdash)$. This fact, together with the claim and $Z \subseteq H / \theta$, implies
\[
H / \theta = Z = H / \theta \cap F / \theta.
\]
Since $\theta$ is compatible both with $H$ and $F$, we conclude that $\bigcap X = H \subseteq F$.
\end{proof}

\section{Meet-prime logics}\label{Sec:meet-prime-logics}

Traditional abstract algebraic logic tends to attribute the status of fundamental concepts both to protoalgebraic and equivalential logics. Unfortunately, this intuition does not match the fact that, when regarded as Leibniz classes, protoalgebraic and equivalential logics happen to be meet-reducible in the Leibniz hierarchy (Theorems \ref{Thm:dec-proto} and \ref{Thm:dec-equivalential}). With an eye towards softening this apparent incoherence, we shall explore a different sense in which a Leibniz class can be considered to capture a primitive or fundamental concept.

\begin{law}
A logic $\vdash$ is \textit{meet-prime} when $\llbracket \vdash \rrbracket$ is meet-prime in $\class{Log}$.
\end{law}

A Leibniz class class can then  also be  considered primitive or fundamental when it is induced by a Leibniz condition whose members are meet-prime logics. In this section we show that this is indeed the case for protoalgebraic, equivalential, and assertional logics.\footnote{For a proof of the fact a range of well-known logics, including all extensions of intuitionistic logic, are meet-irreducible, we refer the reader to \cite{Jan-Mor-19-irr}.}

It is convenient to start with the case of protoalgebraic logics. Recall that for every infinite cardinal $\kappa$, the basic protoalgebraic logic of rank $\kappa$ is denoted by $\vdash_{\class{P}}^{\kappa}$ (Definition \ref{Def:proto-basic}). Our aim is to prove the following:

\begin{Theorem}\label{Thm:basic-proto-prime}
For every infinite cardinal $\kappa > 0$, the logic $\vdash_{\class{P}}^{\kappa}$ is meet-prime.
\end{Theorem}

As a consequence, we obtain the desired result.

\begin{Corollary}
The class of protoalgebraic logics has the form $\mathrm{Log}(\Phi)$ for some Leibniz condition $\Phi$ consisting of meet-prime logics.
\end{Corollary}

\begin{proof}
Immediate from Theorems \ref{Thm:proto-class} and \ref{Thm:basic-proto-prime}.
\end{proof}

The proof of Theorem \ref{Thm:basic-proto-prime} proceeds through a series of technical observations. Given a pair of infinite cardinals $\kappa$ and $\nu$, we let
\[
\{ \langle \A_{j}, F_{j}\rangle \colon j \in J_{\kappa\nu} \}
\]
be the set of $\nu$-generated matrices in $\ModS(\vdash_{\class{P}}^{\kappa})$ up to isomorphism. We can assume without loss of generality  that the various algebras $\A_{j}$ have disjoint universes. Then let $\lambda_{\kappa\nu} \coloneqq \max \{ \omega, \vert \bigcup_{j \in J_{\kappa\nu}} A_{j} \vert\}$ and consider the sets
\begin{align*}
A_{\kappa\nu} & \coloneqq  \{ \top, \bot \} \cup \{ p_{\alpha} \colon \alpha < \lambda_{\kappa\nu}^{+} \}\cup \bigcup_{j \in J_{\kappa\nu}} A_{j}\\
F_{\kappa\nu} & \coloneqq \{ \top \} \cup \{ p_{\alpha} \colon \alpha < \lambda_{\kappa\nu}^{+} \} \cup \bigcup_{j \in J_{\kappa\nu}} F_{j}
\end{align*}
where $\bot, \top, p_{\alpha}$ are new distinct elements. We endow $A_{\kappa\nu}$ with the structure of an $\LL_{\class{P}}^{\kappa}$-algebra $\A_{\kappa\nu}$ stipulating that for every $\alpha < \kappa$, $0 < n \in \omega$, and $a, b, c_{1}, \dots, c_{n} \in A_{\kappa\nu}$,
\begin{displaymath}
a \multimap^{\A_{\kappa\nu}}_{\alpha} b \coloneqq \left\{\begin{array}{@{\,}ll}
a \multimap^{\A_{j}}_{\alpha} b & \text{if $a, b \in A_{j}$ for some $j \in J_{\kappa\nu}$}\\
\top & \text{if $a = b$ and $\{ a, b\} \nsubseteq \bigcup_{j \in J_{\kappa\nu}}A_{j}$}\\
\bot & \text{if $a \ne b$ and $\{ a, b\} \nsubseteq \bigcup_{j \in J_{\kappa\nu}}A_{j}$}\\
\end{array} \right.
\end{displaymath}
and
\begin{displaymath}
\ast_{n \alpha}^{\A_{\kappa\nu}}(c_{1}, \dots, c_{n}) \coloneqq \left\{\begin{array}{@{\,}ll}
\ast_{n \alpha}^{\A_{j}}(c_{1}, \dots, c_{n}) & \text{if $c_{1}, \dots, c_{n} \in A_{j}$ for some $j \in J_{\kappa\nu}$}\\
\bot & \text{otherwise.}\\
\end{array} \right.
\end{displaymath}
Observe that $\A_{\kappa\nu}$ is well defined, since the various $\A_{j}$ have disjoint universes.

\begin{Fact}\label{Fact:proto-1}
For every pair of infinite cardinals $\kappa$ and $\nu$, we have $\vert F_{\kappa\nu} \vert > \vert A_{\kappa\nu} \smallsetminus F_{\kappa\nu} \vert$.
\end{Fact}

\noindent \textit{Proof.}
From the definition of $A_{\kappa\nu}$ and $F_{\kappa\nu}$ it follows that
\[
\pushQED{\qed}\vert F_{\kappa\nu} \vert \geq \lambda^{+}_{\kappa\nu} > \lambda_{\kappa\nu} = \max \{ \omega, \vert \bigcup_{j \in J_{\kappa\nu}} A_{j} \vert\} \geq \vert A_{\kappa\nu} \smallsetminus F_{\kappa\nu} \vert.\qedhere \popQED
\]

\begin{Fact}\label{Fact:proto-2}
For every pair of infinite cardinals $\kappa$ and $\nu$, and $j \in J_{\kappa\nu}$,
\[
\langle \A_{\kappa\nu}, F_{\kappa\nu}\rangle \in \ModS(\vdash_{\class{P}}^{\kappa}) \;\ \text{ and }\;\; \langle \A_{j}, F_{j}\rangle \subseteq \langle \A_{\kappa\nu}, F_{\kappa\nu}\rangle.
\]
\end{Fact}

\begin{proof}
The fact that $\langle \A_{j}, F_{j}\rangle \subseteq \langle \A_{\kappa\nu}, F_{\kappa\nu}\rangle$ is clear. 

To establish that $\langle \A_{\kappa\nu}, F_{\kappa\nu}\rangle \in \ModS(\vdash_{\class{P}}^{\kappa})$, it suffices to prove that the matrix $\langle \A_{\kappa\nu}, F_{\kappa\nu}\rangle$ is a reduced model of $\vdash_{\class{P}}^{\kappa}$. The fact that it is reduced is justified as follows. Consider two distinct elements $a, b \in A_{\kappa\nu}$. We have to prove that $\langle a, b \rangle \notin \leibniz^{\A_{\kappa\nu}}F_{\kappa\nu}$. First we consider the case where there is $j \in J_{\kappa\nu}$ such that $a, b \in A_{j}$. Since $\vdash_{\class{P}}^{\kappa}$ is protoalgebraic and $\langle \A_{j}, F_{j}\rangle \in \ModS(\vdash_{\class{P}}^{\kappa})$, we can apply Theorem \ref{Thm:Protoalgebraic-characterization}(ii), obtaining that the matrix $\langle \A_{j}, F_{j}\rangle$ is reduced. In particular, we can assume without loss of generality that  there is a unary polynomial function $p$ of $\A_{j}$ such that $p(a) \in F_{j}$ and $p(b) \notin F_{j}$ \cite[Prop.\ 2.2(i)]{JaMor19-1}. Since $\langle \A_{j}, F_{j}\rangle \subseteq \langle \A_{\kappa\nu}, F_{\kappa\nu}\rangle$, the map $p$ is also a unary polynomial function of $\A_{\kappa\nu}$ such that $p(a) \in F_{\kappa\nu}$ and $p(b) \notin F_{\kappa\nu}$. As a consequence, $\langle a, b \rangle \notin \leibniz^{\A_{\kappa\nu}}F_{\kappa\nu}$ by \cite[Prop.\ 2.2(i)]{JaMor19-1}. Then we consider the case where the is no $j \in J_{\kappa\nu}$ such that $a, b \in A_{j}$. Choose an arbitrary $\alpha < \kappa$, and consider the unary polynomial function $p(x) \coloneqq x \to_{\alpha}^{\A_{\kappa\nu}}a$ of $\A_{\kappa\nu}$. Observe that
\[
p(b) = a \multimap_{\alpha}^{\A_{\kappa\nu}} b = \bot \notin F_{\kappa\nu}.
\]
On the other hand, we shall see that $p(a) \in F_{\kappa\nu}$. If $a \in A_{j}$ for some $j \in J_{\kappa\nu}$, then we have
\[
p(a) = a \multimap_{\alpha}^{\A_{\kappa\nu}} a \in F_{j} \subseteq F_{\kappa\nu},
\]
since $\langle \A_{j}, F_{j}\rangle \in \Mod(\vdash_{\class{P}}^{\kappa})$ and $\emptyset \vdash_{\class{P}}^{\kappa} x \multimap_{\alpha} x$. Moreover, if $a \notin \bigcup_{j \in J_{\kappa\nu}}A_{j}$, then
\[
p(a) = a \multimap_{\alpha}^{\A_{\kappa\nu}} a = \top \in F_{\kappa\nu}.
\]
Hence we conclude that $p(a) \in F_{\kappa\nu}$. Together with the fact that $p(b) \notin F_{\kappa\nu}$ and \cite[Prop.\ 2.2(i)]{JaMor19-1}, this implies $\langle a, b \rangle \notin \leibniz^{\A_{\kappa\nu}}F_{\kappa\nu}$. We conclude that $\langle \A_{\kappa\nu}, F_{\kappa\nu}\rangle$ is reduced.

It only remains to show that $\langle \A_{\kappa\nu}, F_{\kappa\nu}\rangle \in \Mod(\vdash_{\class{P}}^{\kappa})$, i.e.\ that $\langle \A_{\kappa\nu}, F_{\kappa\nu}\rangle$ is a model of the rules $\emptyset \rhd \Delta_{\kappa}(x, x)$ and $x, \Delta_{\kappa}(x, y) \rhd y$. We detail only the case of $x, \Delta_{\kappa}(x, y) \rhd y$, since the other one is similar. Consider $a, b \in A_{\kappa\nu}$ such that $\{ a \} \cup \Delta^{\A_{\kappa\nu}}(a, b) \subseteq F_{\kappa\nu}$. First we consider the case where $b \in A_{j}$ for some $j \in J_{\kappa\nu}$. Looking at the definition of $\A_{\kappa\nu}$, it is not hard to see that the fact that $b \in A_{j}$ and $\Delta^{\A_{\kappa\nu}}(a, b) \subseteq F_{\kappa\nu}$ implies $a, b \in A_{j}$. In particular, this guarantees that $\Delta^{\A_{\kappa\nu}}(a, b) = \Delta^{\A_{j}}(a, b)$, whence $\Delta^{\A_{j}}(a, b) \subseteq A_{j} \cap F_{\kappa\nu} = F_{j}$. Together with the fact that $\langle \A_{j}, F_{j}\rangle$ is a model of the rule $x, \Delta_{\kappa}(x, y) \rhd y$, this yields $b \in F_{j} \subseteq F_{\kappa\nu}$. Then we consider the case where $b \notin \bigcup_{j \in J_{\kappa\nu}}A_{j}$. Again looking at the definition of $\A_{\kappa\nu}$, it is not difficult to see that the fact that $b \notin \bigcup_{j \in J_{\kappa\nu}}A_{j}$ and $\Delta^{\A_{\kappa\nu}}(a, b) \subseteq F_{\kappa\nu}$ implies $a = b$, whence $b = a \in F_{\kappa\nu}$.
\end{proof}

\begin{Fact}\label{Fact:proto-3}
Let $\kappa$ and $\nu$ be infinite cardinals, $\vdash_{1}$ and $\vdash_{2}$ logics, and $\btau$ an interpretation of $\vdash_{1} \bigotimes \vdash_{2}$ into $\vdash_{\class{P}}^{\kappa}$. Then for every $i = 1, 2$ there is $\langle \B_{i}, G_{i}\rangle \in \ModS(\vdash_{i})$ such that
\[
\langle \A_{\kappa\nu}^{\btau}, F_{\kappa\nu} \rangle \cong \langle \B_{1} \bigotimes \B_{2}, G_{1} \times G_{2}\rangle. 
\]
Moreover, either $\langle \B_{1}, G_{1}\rangle$ or $\langle \B_{2}, G_{2}\rangle$ is trivial.
\end{Fact}

\begin{proof}
Since  $\btau$ is an interpretation of $\vdash_{1} \bigotimes \vdash_{2}$ into $\vdash_{\class{P}}^{\kappa}$, we can apply Fact \ref{Fact:proto-2} obtaining $\langle \A_{\kappa\nu}^{\btau}, F_{\kappa\nu}\rangle \in \ModS(\vdash_{1}\bigotimes \vdash_{2})$. By \cite[Cor.\ 4.14]{JaMor19-1} for every $i = 1, 2$ there is $\langle \B_{i}, G_{i}\rangle \in \ModS(\vdash_{i})$ such that $\langle \A_{\kappa\nu}^{\btau}, F_{\kappa\nu} \rangle \cong \langle \B_{1} \bigotimes \B_{2}, G_{1} \times G_{2}\rangle$. For the sake of simplicity, we assume without loss of generality  that
\begin{equation}\label{Eq:proto-prime-1}
\langle \A_{\kappa\nu}^{\btau}, F_{\kappa\nu} \rangle = \langle \B_{1} \bigotimes \B_{2}, G_{1} \times G_{2}\rangle.
\end{equation}

It only remains to prove that either $\langle \B_{1}, G_{1}\rangle$ or $\langle \B_{2}, G_{2}\rangle$ is trivial. Suppose the contrary, with a view to contradiction. We have 
\begin{align*}
\vert B_{1} \vert + \vert B_{2} \vert &\leq \vert (B_{1} \times (B_{2}\smallsetminus G_{2})) \cup (( B_{1} \smallsetminus G_{1}) \times B_{2}) \vert\\
&= \vert (B_{1} \times B_{2}) \smallsetminus (G_{1} \times G_{2}) \vert\\
&= \vert A_{\kappa\nu} \smallsetminus F_{\kappa\nu} \vert.  
\end{align*}
The first inequality above follows from the fact that, by Lemma \ref{Lem:truth-prime0}, $G_{i} \subsetneq B_{i}$ for every $i = 1, 2$. The second one is obvious, and the third one follows from (\ref{Eq:proto-prime-1}).

Now, recall that $\vert A_{\kappa\nu} \vert \geq \lambda_{\kappa\nu}^{+} \geq \omega$. Thus the set $A_{\kappa\nu} = B_{1} \times B_{2}$ is infinite, whence so is either $B_{1}$ or $B_{2}$. In particular, this implies $\vert B_{1} \vert + \vert B_{2} \vert = \vert B_{1} \times B_{2} \vert = \vert A_{\kappa\nu} \vert$. Together with the above display, this yields
\[
\vert F_{\kappa\nu} \vert \leq \vert A_{\kappa\nu} \vert \leq \vert A_{\kappa\nu} \smallsetminus F_{\kappa\nu} \vert.
\]  
But this is in contradiction with Fact \ref{Fact:proto-1}.
\end{proof}

\begin{proof}[Proof of Theorem \ref{Thm:basic-proto-prime}]
Consider an infinite cardinal $\kappa$. Our aim is to show that the logic $\vdash_{\class{P}}^{\kappa}$ is meet-prime. To this end, consider two logics $\vdash_{1}$ and $\vdash_{2}$ with an interpretation $\btau$ of $\vdash_{1} \bigotimes \vdash_{2}$ into $\vdash_{\class{P}}^{\kappa}$. Then let $\btau_{1}$ be the translation of $\LL_{\vdash_{1}}$ into $\LL_{\class{P}}^{\kappa}$ defined for every $n$-ary $\ast \in \LL_{\vdash_{1}}$ as
\[
\btau_{1}(\ast) \coloneqq \btau(\langle \ast(x_{1}, \dots, x_{n}), x_{1} \rangle). 
\]
The above definition is sound, since the pair $\langle \ast(x_{1}, \dots, x_{n}), x_{1} \rangle$ can be regarded as a basic $n$-ary operation of $\vdash_{1} \bigotimes \vdash_{2}$. Let also $\btau_{2}$ be the translation of $\LL_{\vdash_{2}}$ into $\LL_{\class{P}}^{\kappa}$ defined analogously.

We claim that for every infinite cardinal $\nu$,  there is $i = 1, 2$ such that 
\[
\langle\A_{j}^{\btau_{i}}, F_{j} \rangle \in \ModS(\vdash_{i}) \text{ for every }j \in J_{\kappa\nu}.
\]
To prove this, consider an infinite cardinal $\nu$. By Fact \ref{Fact:proto-3} we can assume without loss of generality  that
\[
\langle \A_{\kappa\nu}^{\btau}, F_{\kappa\nu} \rangle = \langle \A \bigotimes \boldsymbol{1}, F \times \{ 1 \} \rangle
\]
for some $\langle \A, F \rangle \in \ModS(\vdash_{1})$. We shall prove that $\langle\A_{j}^{\btau_{1}}, F_{j} \rangle \in \ModS(\vdash_{1})$ for all $j \in J_{\kappa\nu}$. 
To this end, consider $j \in J_{\kappa\nu}$. By Fact \ref{Fact:proto-2} we obtain $\langle \A_{j}, F_{j}\rangle \subseteq \langle \A_{\kappa\nu}, F_{\kappa\nu}\rangle$, whence
\[
\langle \A_{j}^{\btau}, F_{j}\rangle \subseteq \langle \A_{\kappa\nu}^{\btau}, F_{\kappa\nu}\rangle = \langle\A \bigotimes \boldsymbol{1}, F \times \{ 1 \} \rangle.
\]
Thus there is $\B \subseteq \A$ such that
\begin{equation}\label{Eq:proto-prime-1/2}
\langle \A_{j}^{\btau}, F_{j}\rangle = \langle \B \bigotimes \boldsymbol{1}, (F \cap B) \times \{ 1 \}\rangle.
\end{equation}
Observe that
\begin{equation}\label{Eq:proto-prime-2}
\langle \B, F \cap B \rangle \in \ModS(\vdash_{1}).
\end{equation}
To prove this, observe that $\langle \B, F \cap B \rangle \in \Mod(\vdash_{1})$, since $\langle \B, F \cap B \rangle \subseteq \langle \A, F \rangle$ and $\langle \A, F \rangle \in \Mod(\vdash_{1})$. Hence it only remains to show that $\tarski_{\vdash_{1}}^{\B} F \cap B$ is the identity relation. Consider two distinct elements $a, c \in B$. Since $\langle \A_{j}, F_{j}\rangle \in \ModS(\vdash_{\class{P}}^{\kappa})$ and $\btau$ is an interpretation of $\vdash_{1} \bigotimes \vdash_{2}$ into $\vdash_{\class{P}}^{\kappa}$, we have
\[
\langle \B \bigotimes \boldsymbol{1}, (F \cap B) \times \{ 1 \}\rangle\in \ModS(\vdash_{1} \bigotimes \vdash_{2}).
\]
Since the elements $\langle a, 1 \rangle, \langle c, 1 \rangle \in B \times \{ 1 \}$ are distinct, we can apply \cite[Prop.\ 2.2(ii)]{JaMor19-1} to the above display obtaining without loss of generality  a set $F \cap B \subseteq G \subseteq B$ such that $G \times \{ 1 \}$ is a deductive filter of $\vdash_{1} \bigotimes \vdash_{2}$ on $\B\bigotimes \boldsymbol{1}$, a pair $\langle \varphi(x, y_{1}, \dots, y_{n}), \psi(x, y_{1}, \dots, y_{n})\rangle$ such that $\varphi \in Fm(\vdash_{1})$ and $\psi \in Fm(\vdash_{2})$, and elements $b_{1}, \dots, b_{n} \in B$ such that 
\begin{align*}
\langle \varphi, \psi \rangle^{\B \bigotimes \boldsymbol{1}}(\langle a, 1 \rangle, \langle b_{1}, 1 \rangle, \dots, \langle b_{n}, 1 \rangle) &\in G \times \{ 1 \}\\\
\langle \varphi, \psi \rangle^{\B \bigotimes \boldsymbol{1}}(\langle c, 1 \rangle, \langle b_{1}, 1 \rangle, \dots, \langle b_{n}, 1 \rangle) &\notin G \times \{ 1 \}.
\end{align*}
In particular, we have
\[
\varphi^{\B}(a, b_{1}, \dots, b_{n}) \in G\;\; \text{ and }\;\;  \varphi^{\B}(c, b_{1}, \dots, b_{n}) \notin G.
\]
Observe that $G$ is a deductive filter of $\vdash_{1}$ on $\B$ by Lemma \ref{Lem:min-trick1}. Together with the fact that $F\cap B \subseteq G$ and the above display, this allows us to apply \cite[Prop.\ 2.2(ii)]{JaMor19-1} yielding $\langle a, c \rangle \notin \tarski_{\vdash_{1}}^{\B}(F \cap B)$. This concludes the proof that $\tarski_{\vdash_{1}}^{\B}(F \cap B)$ is the identity relation and establishes (\ref{Eq:proto-prime-2}).

Now, recall from (\ref{Eq:proto-prime-1/2}) that $A_{j} = B \times \{ 1 \}$. Then let $\pi \colon A_{j} \to B$ be the projection on the first coordinate.  The fact that $\pi$ a bijection between $A_{j}$ and $B$ such that $\pi[F_{j}] = F \cap B$ and $F_{j} = \pi^{-1}[F \cap B]$ is a direct consequence of (\ref{Eq:proto-prime-1/2}). Moreover, it is easy to show   that $\pi$ is a homomorphism.  Therefore,  $\pi$ is an isomorphism from $\langle \A^{\btau_{1}}_{j}, F_{j}\rangle$ to $\langle \B, F \cap B \rangle$.
Together with (\ref{Eq:proto-prime-2}) this yields $\langle \A^{\btau_{1}}_{j}, F_{j}\rangle \in \ModS(\vdash_{1})$, establishing the claim.

From the claim it follows that there is $i = 1, 2$ such that for every infinite cardinal $\nu$ there is a cardinal $\mu \geq \nu$ such that $\langle \A_{j}^{\btau_{i}}, F_{j}\rangle \in \ModS(\vdash_{i})$ for all $j \in J_{\kappa\mu}$. Bearing in mind that if $\nu$ and $\mu$ are infinite cardinals such that $\nu \leq \mu$, then $\{\langle \A_{j}, F_{j}\rangle \colon j \in J_{\kappa\nu}\} \subseteq \{\langle \A_{j}, F_{j}\rangle \colon j \in J_{\kappa\mu}\}$, this yields $\langle \A_{j}^{\btau_{i}}, F_{j}\rangle \in \ModS(\vdash_{i})$ for every infinite cardinal $\nu$ and $j \in J_{\kappa\nu}$. Hence we have $\langle \A^{\btau_{i}}, F \rangle \in \ModS(\vdash_{i})$ for every $\langle \A, F \rangle \in \ModS(\vdash_{\class{P}}^{\kappa})$, whence $\btau_{i}$ is a translation of $\vdash_{i}$ into $\vdash_{\class{P}}^{\kappa}$. We conclude that ${\vdash_{i}} \leq {\vdash_{\class{P}}^{\kappa}}$ and, therefore, that $\vdash_{\class{P}}^{\kappa}$ is a meet-prime logic.
\end{proof}

Recall that for every infinite cardinal $\kappa > 0$, the basic equivalential logic of rank $\kappa$ is denoted by $\vdash_{\class{F}}^{\kappa}$ (Definition \ref{Def:equiv-basic}). An argument, similar to the one described above, yields the following conclusion:

\begin{Theorem}\label{Thm:basic-equiv-prime}
For every infinite cardinal $\kappa > 0$, the logic $\vdash_{\class{E}}^{\kappa}$ is meet-prime. As a consequence, the class of equivalential logics has the form $\mathrm{Log}(\Phi)$ for some Leibniz condition $\Phi$ consisting of meet-prime logics.
\end{Theorem}

Recall that the basic assertional logic is denoted by $\vdash_{\class{A}}$ (Definition \ref{Def:asrt-basic}). We have the following:

\begin{Theorem}\label{Thm:basic-equiv-prime}
The logic $\vdash_{\class{A}}$ is meet-prime. As a consequence, the class of assertional logics has the form $\mathrm{Log}(\Phi)$ for a strong Leibniz condition $\Phi$ consisting of a meet-prime logic.
\end{Theorem}

\begin{proof}
In the light to Theorem \ref{Thm:StrongAssertional} it will be enough to show that $\vdash_{\class{A}}$ is meet-prime. To this end, consider two logics $\vdash_{1}$ and $\vdash_{2}$, and an interpretation $\btau$ of $\vdash_{1} \bigotimes \vdash_{2}$ into $\vdash_{\class{A}}$. Moreover, let ${\bf 2}$ be the two-element pointed set. By Proposition \ref{Prop:asrt-class2} we have $\langle {\bf 2}, \{ \top^{\bf 2}_{\ast}\} \rangle \in \ModS(\vdash_{\class{A}})$, whence $\langle {\bf 2}^{\btau}, \{ \top^{\bf 2}_{\ast}\} \rangle \in \ModS(\vdash_{1} \bigotimes \vdash_{2})$. Together with \cite[Cor.\ 4.14]{JaMor19-1}, this implies that for every $i = 1, 2$ there is $\langle \A_{i}, F_{i}\rangle \in \ModS(\vdash_{i})$ such that $\langle {\bf 2}^{\btau}, \{ \top^{\bf 2}_{\ast}\} \rangle = \langle \A_{1} \bigotimes \A_{2}, F_{1} \times F_{2}\rangle$.

Now, observe that $A_{1} \times A_{2}$ is a two-element set. As a consequence, either $A_{1}$ or $A_{2}$ is a singleton. We can assume without loss of generality   that so is $A_{2}$. Together with the fact that $F_{2} \ne \emptyset$, this implies that $\langle \A_{2}, F_{2}\rangle$ is the trivial matrix $\langle {\bf 1}, \{ 1 \} \rangle$. Thus
\[
\langle {\bf 2}^{\btau}, \{ \top^{\bf 2}_{\ast}\} \rangle = \langle \A_{1} \bigotimes {\bf 1}, F_{1} \times \{ 1 \} \rangle.
\]

Consider the translation $\btau_{1}$ of $\LL_{\vdash_{1}}$ into $\LL_{\vdash_{\class{A}}}$ defined as in the proof of Theorem \ref{Thm:basic-proto-prime}. We shall prove that $\btau_{1}$ is an interpretation of $\vdash_{1}$ into $\vdash_{\class{A}}$. Making use of the above display, it is not hard to see that $\langle {\bf 2}^{\btau_{1}}, \{ \top_{\ast}^{\bf 2} \} \rangle \cong \langle \A_{1}, F_{1}\rangle$. Together with the fact that $\langle \A_{1}, F_{1}\rangle \in \ModS(\vdash_{1})$, this yields $\langle {\bf 2}^{\btau_{1}}, \{ \top_{\ast}^{\bf 2} \} \rangle \in \ModS(\vdash_{1})$. By  \cite[Lem.\ 2.3]{JaMor19-1} we conclude
\begin{equation}\label{Eq:asrt-prime-1}
\PSD(\langle {\bf 2}^{\btau_{1}}, \{ \top_{\ast}^{\bf 2} \} \rangle) \subseteq \ModS(\vdash_{1}).
\end{equation}

To prove that $\btau_{1}$ is an interpretation, consider a matrix $\langle \A, F \rangle \in \ModS(\vdash_{\class{A}})$. By Proposition \ref{Prop:asrt-class2} we know that $\A$ is a pointed set and $F = \{ \top_{\ast}^{\A} \}$. In particular, this easily implies $\langle \A, F \rangle \in \PSD(\langle {\bf 2}, \{ \top_{\ast}^{\bf 2} \})$, whence $\langle \A^{\btau_{1}}, F \rangle \in \PSD(\langle {\bf 2}^{\btau_{1}}, \{ \top_{\ast}^{\bf 2} \} \rangle)$. By (\ref{Eq:asrt-prime-1}) this guarantees that $\langle \A^{\btau_{1}}, F \rangle \in \ModS(\vdash_{1})$. Hence we conclude that $\btau_{1}$ is an interpretation of $\vdash_{1}$ into $\vdash_{\class{A}}$, whence ${\vdash_{1}} \leq {\vdash_{\class{A}}}$. This shows that $\vdash_{\class{A}}$ is meet-prime, as desired.
\end{proof}

\paragraph{\bfseries Acknowledgements.}
Thanks are due to James G.\ Raftery for rising the question about whether the theory of the Maltsev and Leibniz hierarchy could be, to some extent, unified. The second author was supported by  the grant CZ.$02$.$2$.$69$/$0$.$0$/$0$.$0$/$17$\_$050$/$0008361$, OPVVV M\v{S}MT, MSCA-IF Lidsk\'{e} zdroje v teoretick\'{e} informatice.

\bibliographystyle{plain}
\end{document}